\documentclass[11pt]{article}

\usepackage{amsmath, amsfonts, amssymb}
\usepackage{lipsum}
\usepackage{graphicx}
\usepackage{epstopdf}
\usepackage{bm}
\usepackage{color}
\usepackage{multirow} 

\newcommand{\bd}[1]{\mathbf{#1}}
\newcommand{\bdg}[1]{\boldsymbol{#1}}
\newcommand{\beq}{\begin{equation}}
\newcommand{\eeq}{\end{equation}}
\newcommand{\pa}{\partial}

\newcommand{\del}{\delta}

\DeclareMathOperator\erf{erf}

\begin{document}

\title{A domain decomposition solution of the Stokes-Darcy system in 3D based on boundary integrals}%
\author{Svetlana Tlupova\thanks{Department of Mathematics, Farmingdale State College, SUNY, Farmingdale, NY 11735, USA tlupovs@farmingdale.edu}}

\date{\today}

\maketitle

\begin{abstract}
A framework is developed for a robust and highly accurate numerical solution of the coupled Stokes-Darcy system in three dimensions. The domain decomposition method is based on a Dirichlet-Neumann type splitting of the interface conditions and solving separate Stokes and Darcy problems iteratively. Second kind boundary integral equations are formulated for each problem. The integral equations use a smoothing of the kernels that achieves high accuracy on the boundary, and a straightforward quadrature to discretize the integrals. Numerical results demonstrate the convergence, accuracy, and dependence on parameter values of the iterative solution for a problem of viscous flow around a porous sphere with a known analytical solution, as well as more general surfaces.
\end{abstract}

{\bf Keywords:} Stokes-Darcy, domain decomposition methods, boundary integral equations, porous sphere.

\section{Introduction}
\label{sec:Intro}

Many important environmental and biological phenomena that occur in our daily lives involve fluids partly flowing freely and partly filtrating through a porous medium. Such processes are of interest in physiology when studying cancer growth and filtration of blood through arterial vessel walls. The free/porous flow plays an important role in numerous industrial applications of air and oil filter design, oil exploration, or chemical reactor simulations. There is also considerable environmental interest in studies of groundwater remediation, geologic CO$_2$ sequestration, and bacterial biofilms.
 
Due to its wide applicability, the modeling of free/porous flow has recently received considerable attention, and mathematical and numerical analyses have been done. The flow is often modeled by the incompressible Stokes equations in the free domain and the Darcy equations in the porous domain, coupled across the interface through suitable conditions. The numerical solutions of these flows have been almost exclusively based on the finite element method~\cite{sun21, li18, discacciati18, chid16, hanspal13, chen11, burman07, correa09, arbogast07, discacciati02, discacciati07, galvis06, guest06, hanspal06, layton03, mardal02, masud02}; a singularity method was used to obtain the solution in~\cite{elasmi01, sekhar00}, and a MAC scheme in~\cite{shiue18}. There are also recently developed partitioned methods for time-­dependent Stokes­-Biot problems, which utilize Robin-­type coupling conditions at the interface~\cite{bukac15, ambart19}.

Despite the significant efforts in recent years to design numerical methods of coupled free/porous flows, these problems remain a challenge today, owing to the multiphysics nature of the system and the large differences in the physical parameters in the governing equations. The severe limitations of the direct solution necessitate the development of techniques that would enable one to split the system and solve the smaller Stokes and Darcy problems separately while exchanging information through the boundary conditions only. These techniques, called domain decomposition methods (DDM), have been used to solve various problems~\cite{quarteroni99}. DDMs are based on dividing the computational domain into several subdomains with or without overlap, and formulating the so-called transmission conditions along the subdomain interfaces. The coupling across the interface is then replaced by an iterative procedure, where the smaller problems are solved independently at each iteration. For the coupled flow problem, similar techniques have been developed and investigated in the context of the finite element method, e.g.~\cite{sun21, cao14, vassilev14, quarteroni99, discacciati02, discacciati07, galvis06}. In particular, we are interested in the Dirichlet-Neumann iterative methods~\cite{quarteroni99, discacciati02}, where the interface conditions are split between the Stokes and Darcy problems. 

In this paper we develop a simple sequential iterative algorithm inspired by the Dirichlet-Neumann methods to decouple the free and porous problems, and formulate highly accurate second-kind boundary integral equations to solve for the fluid quantities. Boundary integral equation methods (BIEM) are extremely powerful in solving many differential equations~\cite{atkinson97, colton98, kress99, hsiao08}, and their importance is well recognized. They have been extensively applied to solve Stokes problems~\cite{pozrikidis92, cortez01, cortez05, hsiao08, greengard96}, as well as porous media flows~\cite{liggett83, rungamornrat06, lough98, pozrikidis03}. In two dimensions, a Robin-Robin DDM was designed for the Stokes-Darcy system and boundary integral equations were used to solve the local problems~\cite{boubendir-tlupova-13}. Non-local operators based on the integral formulations were used in the Robin transmission conditions, leading to a significant improvement in convergence of the DDM. This method is more efficient in overcoming the challenges of the direct solution of~\cite{tlupova09, boubendir09}. However, to the best of our knowledge, no constructive and efficient way of solving the coupled free/porous system in 3D has been developed using BIEM. Attempting a direct solution of the coupled Stokes-Darcy system with BIEM will lead to a formulation that is inevitably of mixed-kind, due to the incompatibility of the Stokes and Darcy differential operators. Once the system is split through the DDM procedure however, suitable formulations can be applied to each problem. The significant advantages of BIEM are the reduced dimensionality as the boundary alone needs to be discretized, high achievable accuracy at points on and off the boundary, the availability of well-conditioned formulations where the efficiency of the solution is maintained with grid refinement, and efficient solutions for problems where the free domain is large or infinite, with boundaries that have complicated geometry, as well as moving boundaries.

To address the singularities that develop in the integrands when solving the integral equations, we use a regularization technique developed in~\cite{beale01, beale04} for Laplace's equation, and recently developed for the Stokes equations in~\cite{tlupova18}. The smoothing can be chosen to have high order on the boundary, so that the accuracy in solving the integral equations generally approaches $O(h^5)$, where $h$ is the grid size in the coordinate planes. This approach is quite simple in that it does not require any special quadrature for evaluating the integrals on the boundary. To discretize the integrals a quadrature rule for closed surfaces from~\cite{wilson, beale16} is applied. This quadrature rule works well for general surfaces without requiring coordinate charts, using projections on coordinate planes instead. As the weights for the quadrature points are found from a partition of unity on the unit sphere, they do not depend on the particular surface. With this approach, the overall formulation is simple to implement, and the data structure needed to describe the boundary is minimal.

The paper is organized as follows. In Section~\ref{sec:ProblemFormulation}, we state the governing equations and suitable interface conditions for a Stokes-Darcy system. In Section~\ref{sec:DDM}, we outline a Dirichlet-Neumann type iterative procedure to split and sequentially solve the Darcy and Stokes problems. A convergence analysis of this algorithm for a spherical geometry based on spherical harmonics representations is conducted in Section~\ref{sec:Analysis}. In Section~\ref{sec:BIEM}, a boundary integral formulation is developed for both Darcy and Stokes equations, including a description of the treatment of singularities in the integrands in Section~\ref{sec:Regularization} and the quadrature in Section~\ref{sec:Quadrature}. The complete algorithm is presented in Section~\ref{sec:Algorithm} and the discrete iteration operator is computed in Section~\ref{sec:DiscIterOper}. Numerical results are presented in Section~\ref{sec:NumericalResults}, where we use a known analytical solution for a viscous flow around a porous sphere as a benchmark problem, as well as test the algorithm on other geometries.


\section{Problem formulation}
\label{sec:ProblemFormulation}

\begin{figure}[h]
\begin{centering}
\scalebox{0.75}{\includegraphics{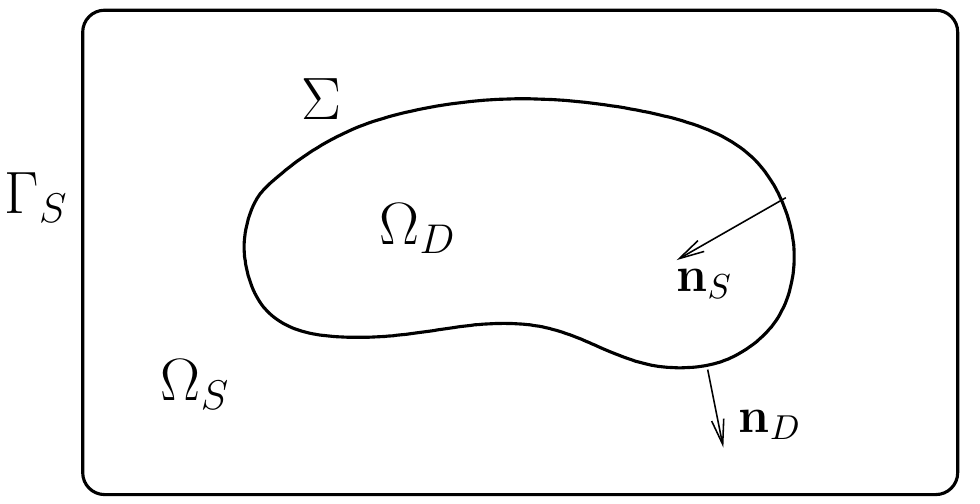}} 
\caption{The domain with the free fluid region  $\Omega_S$ and the porous medium region $\Omega_D$, with the common interface $\Sigma$.}
\label{Schematic}
\end{centering}
\end{figure}

To describe the model equations, we denote the free flow quantities and the porous domain quantities by subscripts $S$ and $D$, respectively. Figure~\ref{Schematic} shows an example two-dimensional setting, where a free fluid flow in a bounded region $\Omega_S$ and a flow in a porous region $\Omega_D$ are coupled across the common interface $\Sigma$.
In the region of free fluid flow $\Omega_S$, we assume the low Reynolds number regime modeled by the incompressible Stokes equations,
\beq
\label{StokesEqs1}
\textrm{In} \ \Omega_{S} :
\left \{ \begin{array}{l}
        \nabla p_S - \mu \Delta\bd{u}_S = 0, \\[2pt]
          \nabla \cdot \bd{u}_S = 0,
\end{array}\right.
\eeq
where $\mu$ is the fluid viscosity, $p_S$ and $\bd{u}_S$ are the fluid pressure and velocity, respectively. Darcy's equations govern the flow in the porous medium,
\beq
 \label{DarcyEqs1}
 \textrm{In} \ \Omega_{D} :
 \left \{ \begin{array}{l}
    \nabla p_D + \mu K^{-1}\bd{u}_D = 0, \\[2pt]
   \nabla \cdot \bd{u}_D = 0.
   \end{array}\right.
\eeq
Here $\bd{u}_D$ and $p_D$ are respectively  the (averaged) fluid velocity and the hydrostatic pressure, and $K$ is a symmetric and positive definite permeability tensor that depends on the microstructure of the porous medium. The medium is called homogeneous if the permeability $K(\bd{x})$ is constant throughout the domain, and isotropic when permeability is a scalar, $K = \kappa I$. In this work, we treat porous media that are homogeneous and isotropic.

The interface between these two systems of partial differential equations serves as a replacement of the boundary layer through which the velocity changes rapidly. A lot of attention has been devoted to the formulation and analysis of appropriate coupling conditions~\cite{beavers67,saffman71, jager96, jager00, ochoa-tapia95_1, ochoa-tapia95_2, sahraoui92}. The incompressibility condition leads to continuity of normal velocity components,
\beq
    \label{Vel_contin}
    \bd{u}_S \cdot \bd{n}_S = -\bd{u}_D \cdot \bd{n}_D,
\eeq
where $\bd{n}_S$ and $\bd{n}_D$ are the unit normal vectors that point out of the regions $\Omega_D$ and $\Omega_S$ respectively, so that $\bd{n}_D = -\bd{n}_S$ on $\Sigma$ (see Fig.~\ref{Schematic}). In addition, a suitable condition on the tangential velocity is the slip condition of Beavers-Joseph-Saffman~\cite{beavers67, saffman71},
\beq
	\label{BeaversJoseph}
	\bdg{\tau}_S\cdot\bdg{\sigma}_S\cdot\bd{n}_S = \frac{\gamma}{\sqrt{\kappa}}\bd{u}_S \cdot \bdg{\tau}_S,
\eeq
where $\bdg\sigma_S = -p_S I+\mu[\nabla\bd{u}_S + (\nabla\bd{u}_S)^T]$ is the Stokes stress tensor, $\bdg{\tau}_S$ is the tangent vector to the interface, and $\gamma$ is a dimensionless quantity that depends on the structure of the porous material. Well-posedness of the mathematical formulation with this condition has been analyzed for steady flows~\cite{burman07, discacciati02, galvis06, layton03} as well as unsteady flows~\cite{cao09}. The final condition represents the balance of normal forces across the interface,
\beq
    \label{Press_contin}
    \bd{n}_S\cdot\bdg{\sigma}_S\cdot\bd{n}_S = -p_D,
\eeq
which allows the pressure to be discontinuous across the interface.

A fundamental difficulty in the development of robust computational methods for the free/porous flow system~\eqref{StokesEqs1}-\eqref{Press_contin} is due to the fact that the PDEs that govern the fluid flow in the free region and the flow in the porous medium have different properties. The differential operators in the two subdomains are incompatible, inhibiting the simultaneous solution of the equations by a direct method. In addition, the equations include parameters, the viscosity $\mu$ and permeability $\kappa$, that could differ by several orders of magnitude in practice.


\section{The Dirichlet-Neumann iterative method}
\label{sec:DDM}

The severe limitations of the direct solution necessitate the development of techniques that would enable one to split the system and solve the smaller Stokes and Darcy problems separately while exchanging information through the boundary conditions only. We develop a numerical scheme based on the sequential Dirichlet-Neumann (D-N) iterative method~\cite{quarteroni99, discacciati02}. The outline of the algorithm is as follows. First, we define $q = -\bd{u}_S\cdot \bd{n}_S$ as the iteration variable on the interface $\Sigma$, and choose an initial guess $q^{(0)} = 0$. Then, for $k=1,2,3,...$ until convergence, perform the following steps:
\begin{enumerate}
\item Solve the Darcy problem with the interface condition~\eqref{Vel_contin},
	\beq
	\label{DarcyProblem1}
 	\left \{ \begin{array}{ll}
    	\nabla p^{(k)}_D + \mu \kappa^{-1}\bd{u}^{(k)}_D = 0 & \textrm{in} \ \Omega_{D}, \\[2pt]
   	\nabla \cdot \bd{u}^{(k)}_D = 0 & \textrm{in} \ \Omega_{D},\\[2pt]
	\bd{u}^{(k)}_D \cdot \bd{n}_D = q^{(k-1)} & \textrm{on} \ \Sigma.
   	\end{array}\right.
	\eeq
	The interface condition is a Neumann condition for $p_D$. 
\item Using $p^{(k)}_D$ from Step 1, solve the Stokes problem with the interface conditions~\eqref{BeaversJoseph}-\eqref{Press_contin},
	\beq
	\label{StokesProblem1}
 	\left \{ \begin{array}{ll}
    	\nabla p^{(k)}_S - \mu \Delta\bd{u}^{(k)}_S = \bd{0} & \textrm{in} \ \Omega_{S}, \\[2pt]
   	\nabla \cdot \bd{u}^{(k)}_S = 0 & \textrm{in} \ \Omega_{S},\\[2pt]
	\bdg{\tau}_S\cdot\bdg{\sigma}^{(k)}_S\cdot\bd{n}_S = \frac{\gamma}{\sqrt{\kappa}}\bd{u}^{(k-1)}_S\cdot\bdg{\tau}_S & \textrm{on} \ \Sigma,\\[2pt]
	\bd{n}_S\cdot\bdg{\sigma}^{(k)}_S\cdot\bd{n}_S = -p^{(k)}_D & \textrm{on} \ \Sigma.
   	\end{array}\right.
	\eeq
\item Update the iteration variable
	\beq
	\label{q_update}
	q^{(k)} = (1-\theta) \, q^{(k-1)} - \theta \, \bd{u}^{(k)}_S\cdot \bd{n}_S
	\eeq 
	using $\bd{u}^{(k)}_S$ from Step 2, where $\theta \in (0,1)$ is a relaxation parameter.
\end{enumerate}
The convergence criterion is taken as the relative error between $q^{(k+1)}$ and $q^{(k)}$, which we call the DDM residual, falling below a prescribed tolerance. The solutions of the now separate Darcy \eqref{DarcyProblem1} and Stokes \eqref{StokesProblem1} problems are based on the boundary integral formulations which we describe in Section~\ref{sec:BIEM}.


\section{Convergence analysis}
\label{sec:Analysis}

We analyze the convergence properties of the Dirichlet-Neumann iterative scheme using the spherical harmonics representation of both problems, assuming as the interface a sphere of radius $R$ centered at the origin, with the Darcy problem inside and the Stokes problem outside the sphere.

\subsection{Darcy problem inside a sphere}

We consider the problem
\beq
\label{CA1}
 	\left \{ \begin{array}{ll}
    	\tilde{\kappa}\, \nabla p + \bd{u} = 0 & \textrm{in} \ \Omega, \\[2pt]
   	\nabla \cdot \bd{u} = 0 & \textrm{in} \ \Omega,\\[2pt]
	\bd{u} \cdot \bd{n} = q & \textrm{on} \ \Sigma,
   	\end{array}\right.
\eeq
where $\Omega$ is the sphere of radius $R$ centered at the origin, $\Sigma$ is its boundary, $\bd{n}$ is the outward unit normal vector, and $\tilde{\kappa} = \kappa/\mu$.

In spherical coordinates $(r,\theta,\phi)$, the general solution to Laplace's equation $\Delta p=0$ in $\Omega$ can be written as 
\beq
\label{CA2}
	p = \sum_{n=0}^\infty p_n,
\eeq
with
\beq
\label{CA3}
	p_n = \sum_{m=-n}^n a_n^m\, r^n\, Y_n^m(\theta,\phi),
\eeq
where $p_n$ is a solid harmonic of order $n$, $Y_n^m$ is the spherical harmonic function of degree $n$ and order $m$, and $a_n^m\in\mathbb{C}$ are constants.

With $\bd{n}=\bd{x}/r$ and 
\beq
\label{CA4}
	\nabla p = \sum_{n=0}^\infty \sum_{m=-n}^n \left[ \frac{\pa}{\pa r} (a_n^m\, r^n)\, \frac{\bd{x}}{r} \, Y_n^m(\theta,\phi) + a_n^m\, r^n\, \nabla Y_n^m(\theta,\phi)\right],
\eeq
for the normal velocity we get
\beq
\label{CA5}
	\bd{u}\cdot\bd{n} = -\tilde{\kappa}\, \nabla p \cdot \bd{n} = -\tilde{\kappa} \sum_{n=0}^\infty \sum_{m=-n}^n \left[ n\, a_n^m\, r^{n-1}\, Y_n^m(\theta,\phi) \right].
\eeq

Using a similar representation for $q$ on $\Sigma$, 
\beq
\label{CA6}
	q = \sum_{n=0}^\infty \sum_{m=-n}^n q_n^m\, Y_n^m(\theta,\phi),
\eeq
and imposing the Neumann boundary condition in~\eqref{CA1} on $r=R$, we find for the coefficients $a_n^m$,
\beq
\label{CA7}
	-\tilde{\kappa}\, n\, a_n^m\, R^{n-1} = q_n^m.
\eeq


\subsection{Stokes problem outside a sphere}

We now consider the solution of the Stokes problem with a no-slip boundary condition for simplicity,
\beq
\label{CA8}
 	\left \{ \begin{array}{ll}
    	\nabla p - \mu\, \Delta\bd{u} = 0 & \textrm{in} \ \Omega^c, \\[2pt]
   	\nabla \cdot \bd{u} = 0 & \textrm{in} \ \Omega^c,\\[2pt]
   	\bd{u}\cdot\bdg{\tau} = 0 & \textrm{on} \ \Sigma,\\[2pt]
	\bd{n}\cdot \bdg{\sigma} \cdot \bd{n} = -p_D & \textrm{on} \ \Sigma,
   	\end{array}\right.
\eeq
where $p_D$ is the Darcy pressure on $\Sigma$ and $\Omega^c$ is the exterior of the sphere of radius $R$.

The solution to the Stokes equations in spherical harmonics is referred to as Lamb's general solution~\cite{lamb32}. Since the Stokes pressure satisfies Laplace's equation, it can be represented similar to~\eqref{CA2}-\eqref{CA3}, but for the exterior solution we keep only the negative harmonics,
\beq
\label{CA9}
	p = \sum_{n=0}^\infty p_{-n-1},
\eeq
with
\beq
\label{CA10}
	p_{-n-1} = \sum_{m=-n}^n b_n^m\, r^{-n-1}\, Y_n^m(\theta,\phi).
\eeq
The disturbance velocity field is represented by (see~\cite{kim-karrila05})
\beq
\label{CA11}
	\bd{u} = \sum_{n=1}^\infty \left[ -\frac{(n-2)\, r^2\, \nabla p_{-n-1}}{2\mu\, n (2n-1)} + \frac{(n+1)\, \bd{x}\, p_{-n-1}}{\mu\, n (2n-1)} + \nabla \phi_{-n-1} + \nabla \times (\bd{x}\, \chi_{-n-1}) \right],
\eeq
where
\beq
\label{CA12}
	\phi_{-n-1} = \sum_{m=-n}^n c_n^m\, r^{-n-1}\, Y_n^m(\theta,\phi),
\eeq
\beq
\label{CA13}
	\chi_{-n-1} = \sum_{m=-n}^n d_n^m\, r^{-n-1}\, Y_n^m(\theta,\phi).
\eeq
The first two terms in~\eqref{CA11} correspond to a particular solution, while the rest corresponds to the homogeneous solution, constructed from a potential (the $\nabla\phi$ term) and a toroidal field (the $\nabla \times (\bd{x}\, \chi)$ term)~\cite{kim-karrila05}.

The stress vector acting across the sphere can be expressed in general as (see~\cite{happel-brenner83} equation (3-2.37) or~\cite{kim-karrila05} exercise 4.6, p.104),
\begin{eqnarray}
\label{CA14}
	\bdg{\sigma}_n = \bdg{\sigma}\cdot \bd{n} = \frac{1}{r} \sum_{n=-\infty}^\infty \left[ \frac{n(n+2)\, r^2\, \nabla p_n}{(n+1) (2n+3)} - \frac{(2n^2+4n+3)\, \bd{x}\, p_n}{(n+1) (2n+3)} \right] \nonumber\\ 
	+ \frac{\mu}{r} \sum_{n=-\infty}^\infty \left[ 2(n-1)\nabla \phi_n + (n-1)\nabla \times (\bd{x}\, \chi_n) \right].
\end{eqnarray}
For the exterior problem, we again only keep the negative harmonics and replace $n$ in~\eqref{CA14} with $-n-1$, to obtain
\begin{eqnarray}
\label{CA15}
	\bdg{\sigma}_n = \bdg{\sigma}\cdot \bd{n} = \frac{1}{r} \sum_{n=1}^\infty \left[ \frac{(n-1)(n+1)\, r^2\, \nabla p_{-n-1}}{n (2n-1)} - \frac{(2n^2+1)\, \bd{x}\, p_{-n-1}}{n (2n-1)} \right] \nonumber\\ 
	- \frac{\mu}{r}\, \sum_{n=1}^\infty \left[ 2(n+2)\nabla \phi_{-n-1} + (n+2)\nabla \times (\bd{x}\, \chi_{-n-1}) \right].
\end{eqnarray}

We now turn to the boundary conditions. The no-slip condition on $r=R$ implies (see~\cite{kim-karrila05} p.89) that $\chi$ is not needed and
\beq
\label{CA16}
	\sum_{n=1}^\infty \left[ -\frac{n(n+1)\, R}{2\mu\, (2n-1)}\, p_{-n-1} |_{r=R} + \frac{(n+1)(n+2)}{R} \phi_{-n-1} |_{r=R} \right] = 0,
\eeq
which leads to
\beq
\label{CA17}
	\phi_{-n-1} |_{r=R} = \frac{n\, R^2}{2\mu\, (n+2)(2n-1)} p_{-n-1} |_{r=R}.
\eeq
Note that for the radial component of velocity we have from~\eqref{CA11},
\beq
\label{CA18}
	u_r = \sum_{n=1}^\infty \left[ -\frac{(n-2)\, r^2}{2\mu\, n (2n-1)}\, \frac{\pa p_{-n-1}}{\pa r} + \frac{(n+1)\, r\, p_{-n-1}}{\mu\, n (2n-1)} + \frac{\pa \phi_{-n-1}}{\pa r} \right].
\eeq
Using Euler's theorem for homogeneous polynomials,
\beq
\label{CA19}
	r\, \frac{\pa h_n}{\pa r} = n\, h_n,
\eeq
where $h_n$ is any solid harmonic of order $n$, the radial velocity in~\eqref{CA18} simplifies to
\beq
\label{CA20}
	u_r = \sum_{n=1}^\infty \left[ \frac{(n+1)\, r\, p_{-n-1}}{2\mu\, (2n-1)} - \frac{(n+1)}{r} \phi_{-n-1} \right].
\eeq

Similarly, for the radial component of the stress vector $\bdg{\sigma}_n$,
\begin{eqnarray}
\label{CA21}
	\sigma_{nr} = \frac{1}{r} \sum_{n=1}^\infty \left[ \frac{(n-1)(n+1)\, r^2}{n (2n-1)}\, \frac{\pa p_{-n-1}}{\pa r} - \frac{(2n^2+1)\, r\, p_{-n-1}}{n (2n-1)} \right] \nonumber\\ 
	- \frac{\mu}{r}\, \sum_{n=1}^\infty \left[ 2(n+2) \frac{\pa \phi_{-n-1}}{\pa r} \right],
\end{eqnarray}
which simplifies again using~\eqref{CA19},
\beq
\label{CA22}
	\sigma_{nr} = \sum_{n=1}^\infty \left[ -\frac{(n^2+3n-1)}{(2n-1)} p_{-n-1} + \frac{2\mu}{r^2} (n+1)(n+2) \phi_{-n-1} \right].
\eeq
On the boundary $r=R$ using~\eqref{CA17}, this becomes
\beq
\label{CA23}
	\sigma_{nr}|_{r=R} = \sum_{n=1}^\infty \left[ - p_{-n-1} \right]_{r=R}, 
\eeq
and therefore enforcing the boundary condition on the normal stress
\beq
\label{CA24}
	\bdg{\sigma}_n\cdot\bd{n} = \sigma_{nr}|_{r=R} = -p_D 
\eeq
leads to
\beq
\label{CA25}
	\sum_{n=1}^\infty \sum_{m=-n}^n b_n^m\, R^{-n-1}\, Y_n^m(\theta,\phi) = \sum_{n=1}^\infty \sum_{m=-n}^n a_n^m\, R^n\, Y_n^m(\theta,\phi).
\eeq
Using~\eqref{CA7}, we write
\beq
\label{CA26}
	b_n^m = -\frac{R^{n+2}}{\tilde{\kappa}\, n} q_n^m.
\eeq

\subsection{Spectrum of the Dirichlet-Neumann procedure}

We now diagonalize the iteration operator of the Dirichlet-Neumann method defined on the sphere. To derive the spectrum of the iteration operator, note that the iteration variable $q$ is updated at each iteration according to~\eqref{q_update}, so we need the normal component of Stokes velocity on the sphere. We evaluate the radial velocity~\eqref{CA20} on the boundary $r=R$ and use~\eqref{CA17} to get for the normal Stokes velocity,
\beq
\label{CA27}
	u_r |_{r=R} = \sum_{n=1}^\infty \left[ \frac{(n+1)\, R}{\mu\, (n+2)(2n-1)} p_{-n-1}|_{r=R} \right].
\eeq
Using~\eqref{CA26}, this becomes
\beq
\label{CA28}
	u_r |_{r=R} = -\frac{R^2}{\kappa} \sum_{n=1}^\infty \left[ \frac{(n+1)}{n(n+2)(2n-1)} \sum_{m=-n}^n q_n^m\,  Y_n^m(\theta,\phi) \right].
\eeq
For the DDM update we get
\beq
\label{CA29}
	q^{(k)} = (1-\theta) q^{(k-1)} + \theta u_r^{(k)} |_{r=R},
\eeq
which we write in the form of the diagonalized iteration operator $q^{(k)} = \mathcal{A} q^{(k-1)}$ as
\beq
\label{CA30}
	q_n^{m(k)} = \mathcal{A}_n q_n^{m(k-1)},
\eeq
where
\beq
\label{CA31}
	\mathcal{A}_n = (1-\theta) - \theta\frac{R^2}{\kappa} \frac{(n+1)}{n(n+2)(2n-1)}, \qquad n=1,2,3,...
\eeq

\begin{figure}[htb]
\begin{centering}
	\scalebox{0.27}{\includegraphics{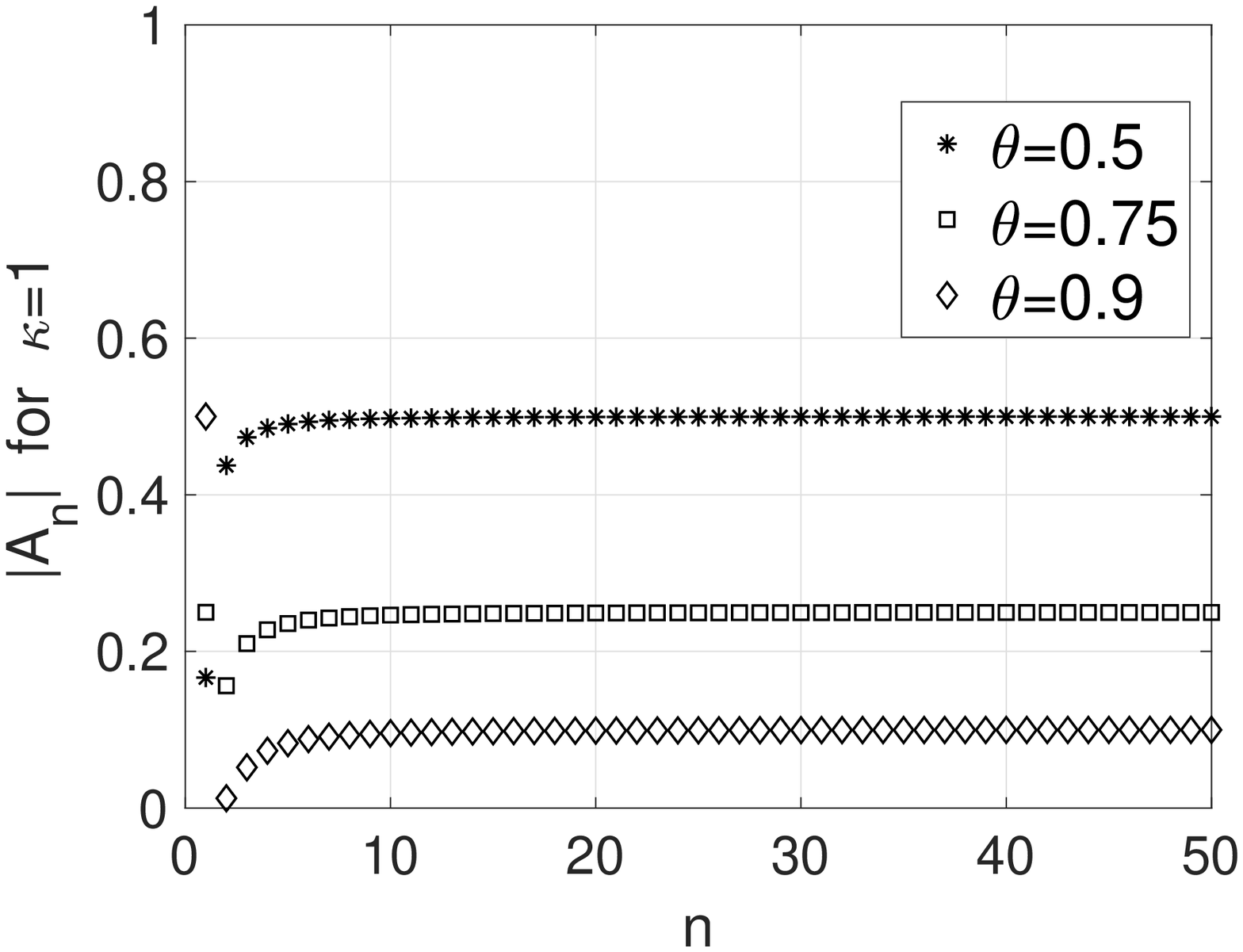}}
	\scalebox{0.27}{\includegraphics{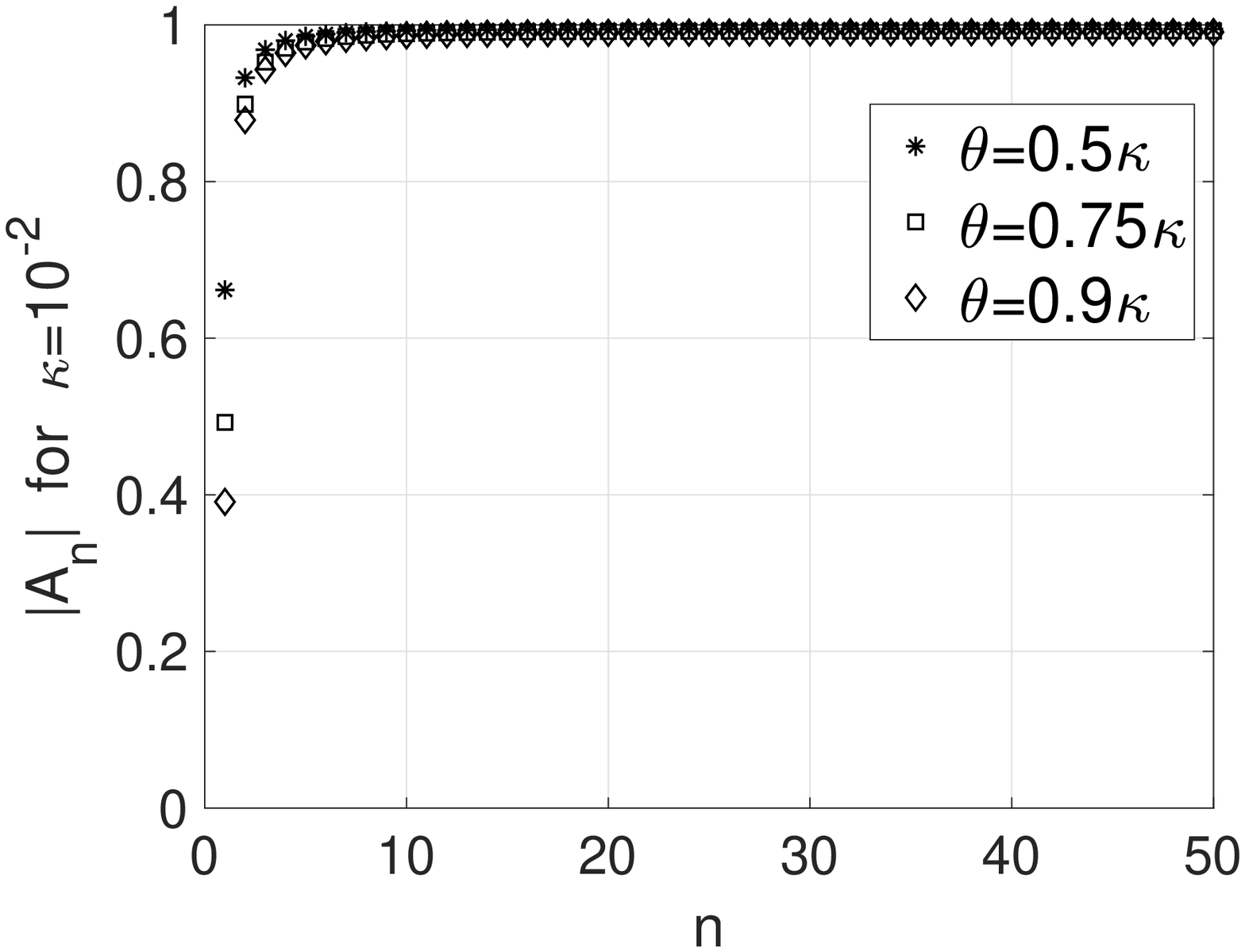}}
	\scalebox{0.27}{\includegraphics{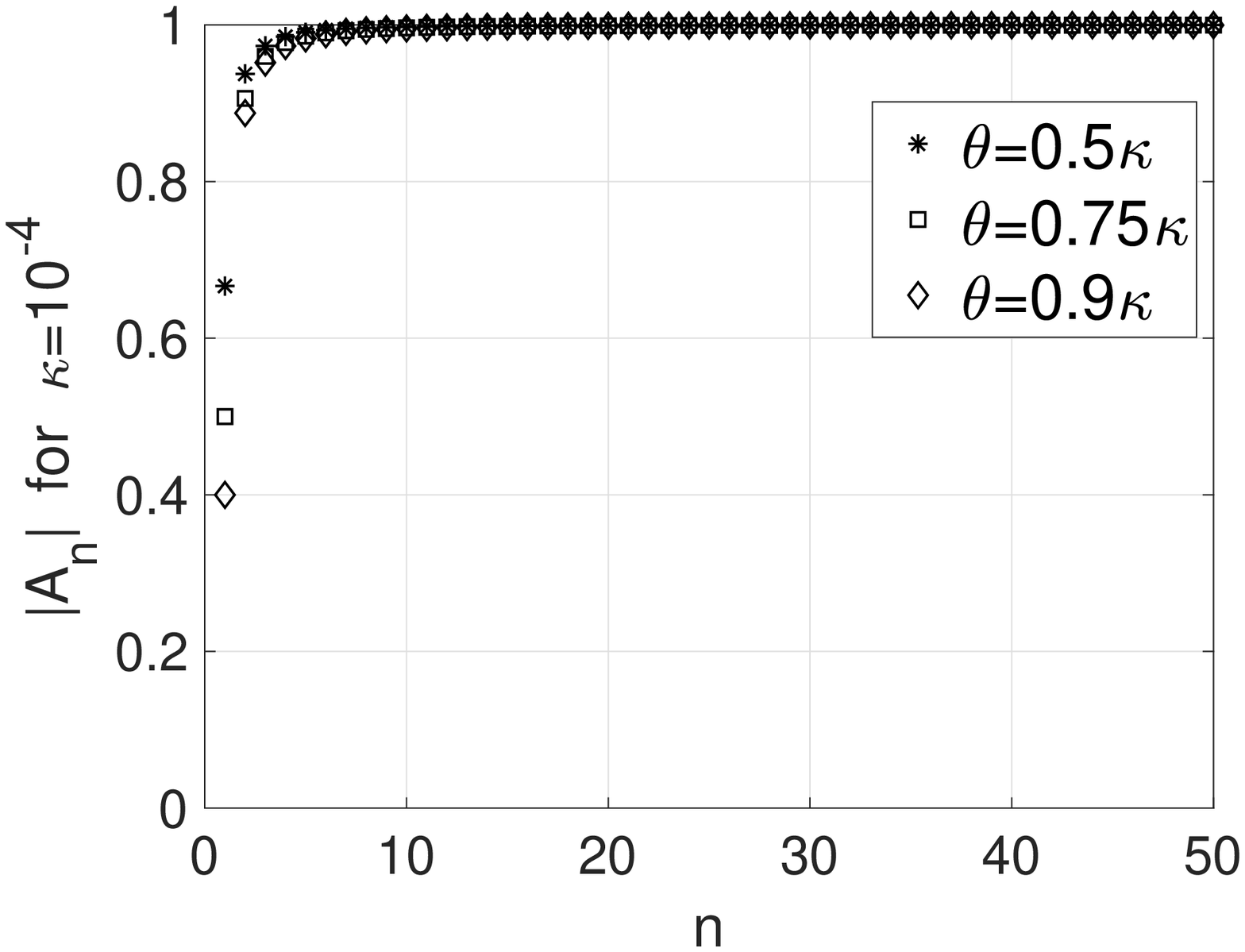}}
	\caption{Harmonic coefficients $\mathcal{A}_n$ from~\eqref{CA31} of the iteration operator on the unit sphere, for permeability values $\kappa=1, 10^{-2}, 10^{-4}$.}
\label{Sphere_modes}
\end{centering}
\end{figure}

Figure~\ref{Sphere_modes} shows these coefficients for $n=1,...,50$ for different values of permeability $\kappa$ as well as different $\theta$. We can observe that for $\kappa=1$, a larger value of $\theta$ leads to smaller coefficients $\mathcal{A}_n$ and therefore the iterations are expected to converge faster, as we will show in the numerical tests. For smaller values of $\kappa$ however, the coefficients $\mathcal{A}_n$ get close to 1 as $n$ increases, independent of the value of $\theta$. This can hinder the convergence rate of the algorithm when using a successive approximation technique as in~\eqref{q_update}. We will show in our numerical tests that this indeed is what happens, and switch to using a more efficient method, and one that does not require all eigenvalues to be less than 1, such as GMRES to help mitigate this issue. Interesting approaches that will be investigated in future work are (i) different decoupling, e.g. Robin-Robin transmission conditions, and (ii) using in place of $\theta$ non-local operators based on the boundary integral formulation to speed up the convergence rates. Similar techniques were shown to dramatically speed up the iterations in two dimensions~\cite{boubendir-tlupova-13}.


\section{Boundary integral formulations}
\label{sec:BIEM}

There are two main approaches to formulating boundary integral equations, usually referred to as the 'direct' and 'indirect' approaches. The direct approach is based on the Green's representation formula and is the approach we take in this work. While the analytic properties of the integral operators are the same with both approaches, our primary interest in the direct approach is that the resulting equations are for the physically relevant Dirichlet or Neumann data.


\subsection{Darcy problem}

The Green's representation formula,
\beq
	\label{Greens}
	p_D(\bd{y}) = - \int_{\pa\Omega_D} \frac{\pa p_D(\bd{x})}{\pa n(\bd{x})} \,G(\bd{y,x})\, dS(\bd{x}) + \int_{\pa\Omega_D} p_D(\bd{x}) \, \frac{\pa G(\bd{y,x})}{\pa n(\bd{x})} \, dS(\bd{x}) , \qquad \bd{y}\in\Omega_D,
\eeq
represents any harmonic function $p_D$ in terms of its Cauchy data, that is, its boundary values and its normal derivative on the boundary. The first integral in~\eqref{Greens} is the single-layer potential and the second integral is the double-layer potential, where the kernels are $G(\bd{y,x}) = -1/4\pi |\bd{y} - \bd{x}|$ and its derivative in the outward normal $\frac{\pa G(\bd{y,x})}{\pa n(\bd{x})} = \nabla_\bd{x} G(\bd{y},\bd{x})\cdot \bd{n}(\bd{x}) = -(\bd{y}-\bd{x})\cdot \bd{n}(\bd{x})/4\pi |\bd{y}-\bd{x}|^3$. The double layer potential undergoes a discontinuity across the boundary, and the jump relations 
\beq
	\label{Jump}
	\int_{\pa\Omega_D} \frac{\pa G(\bd{y,x})}{\pa n(\bd{x})} \, dS(\bd{x}) 
	= \left\{ \begin{array}{ll} 
	1, & \bd{y}\in\Omega_D\\
	1/2, & \bd{y}\in\pa\Omega_D\\
	0, & \bd{y}\in\mathbb{R}^3 \setminus \bar{\Omega}_D
	\end{array}\right.
\eeq
are used to extend~\eqref{Greens} to the boundary,
\beq
	\label{Greens_bdry}
	\frac12 p_D(\bd{y}) = - \int_{\pa\Omega_D} \frac{\pa p_D(\bd{x})}{\pa n(\bd{x})} \,G(\bd{y,x})\, dS(\bd{x}) + \int_{\pa\Omega_D} p_D(\bd{x}) \, \frac{\pa G(\bd{y,x})}{\pa n(\bd{x})} \, dS(\bd{x}) , \qquad \bd{y}\in\pa\Omega_D.
\eeq
We write~\eqref{Greens_bdry} as a second-kind integral equation assuming known Neumann data $\pa p_D/\pa n$ on the boundary,
\beq
	\label{Darcy_BIEM}
	\frac12 p_D = H p_D + b_D, \qquad \qquad b_D = -L \frac{\pa p_D}{\pa n},
\eeq
where by $L$ and $H$ we denote the integral operators corresponding to the single- and double-layer potentials, respectively. The integral equation~\eqref{Darcy_BIEM} can be solved by successive approximations, otherwise known as a Neumann's iteration scheme (see, for example,~\cite{kress99} or~\cite{hsiao08}, Sec. 5.6.7),
\beq
	\label{Darcy_BIEM_iter}
	p_D^{(n+1)} = \frac12 p_D^{(n)} + H p_D^{(n)} + b_D, \qquad n=0,1,2,...
\eeq


\subsection{Stokes problem}

We follow a similar 'direct' approach for the Stokes problem and use the following integral representation for the velocity,
\beq
	\label{Stokes}
	\bd{u}_S(\bd{y}) = \bd{u}^\infty -\frac{1}{8\pi}\int_{\pa\Omega_S} S(\bd{y,x})\, \bd{f}_S(\bd{x}) \, dS(\bd{x}) + \frac{1}{8\pi}\int_{\pa\Omega_S} \bd{u}_S(\bd{x}) \cdot T(\bd{y,x}) \cdot \bd{n}_S(\bd{x}) \, dS(\bd{x}) , \quad \bd{y}\in\Omega_S,
\eeq
where $\bd{n}_S$ is the unit outward normal vector to the Stokes domain boundary, $\bd{f}_S = \bdg{\sigma}_S\cdot\bd{n}_S$ is the surface force, or traction, on the boundary $\pa\Omega_S$, $\bdg{\sigma}_S$ is the stress tensor, and
\begin{equation}
	S_{ij}(\bd{y,x}) = \frac{\del_{ij}}{|\bd{\hat{y}}|} + \frac{\hat{y}_i \ \hat{y}_j}{|\bd{\hat{y}}|^3}, \qquad \qquad
	T_{ijk} (\bd{y,x}) = -6 \frac{ \hat{y}_i\hat{y}_j\hat{y}_k}{|\bd{\hat{y}}|^5},
\end{equation}
are the Stokeslet and stresslet fundamental solutions of the Stokes equations, where $\bd{\hat{y}} = \bd{y}- \bd{x}$, $\del_{ij}$ is the Kronecker delta, and $i,j,k = 1,2,3$ are Cartesian coordinates. A background flow $\bd{u}^\infty$ is added for external flows. 

The representation~\eqref{Stokes} is based on the Lorentz reciprocal theorem, and again, the first integral is called the single-layer potential and the second integral is the double-layer potential of Stokes flow. The double-layer undergoes a discontinuity across the boundary much like~\eqref{Jump}, and we again use jump conditions to continuously extend the representation~\eqref{Stokes} to the boundary:
\beq
	\label{Stokes_bdry}
	\frac12 \bd{u}_S(\bd{y}) = \bd{u}^\infty -\frac{1}{8\pi}\int_{\pa\Omega_S} S(\bd{y,x})\, \bd{f}_S(\bd{x}) \, dS(\bd{x}) + \frac{1}{8\pi}\int_{\pa\Omega_S} \bd{u}_S(\bd{x}) \cdot T(\bd{y,x}) \cdot \bd{n}_S(\bd{x}) \, dS(\bd{x}) , \quad \bd{y}\in\pa\Omega_S.
\eeq
We write~\eqref{Stokes_bdry} as 
\beq
	\label{Stokes_BIEM}
	\frac12 \bd{u}_S = K \bd{u}_S + b_S, \qquad \qquad b_S = -M \bd{f}_S + \bd{u}^\infty,
\eeq
where by $M$ and $K$ we denote the integral operators corresponding to the Stokes single- and double-layer potentials, respectively. Assuming a known surface force, equation~\eqref{Stokes_BIEM} is a second-kind integral equation for $\bd{u}_S$, and we solve it by successive approximations
\beq
	\label{Stokes_BIEM_iter}
	\bd{u}_S^{(n+1)} = \frac12 \bd{u}_S^{(n)} + K \bd{u}_S^{(n)} + b_S, \qquad n=0,1,2,...
\eeq

It is a well-known challenge in the numerical solution of the integral equations~\eqref{Darcy_BIEM_iter} and~\eqref{Stokes_BIEM_iter} that singularities develop in the integrands $H, L, K, M$, and we address this issue by using a high accuracy regularization technique which we describe next.


\subsection{High-order regularization of kernels on the surface}
\label{sec:Regularization}

The integrands in~\eqref{Greens_bdry} and~\eqref{Stokes_bdry} exhibit singularities for $\bd{y}=\bd{x}$. To address this issue, we first use subtraction in the double layer integrals to reduce the singularities,
\begin{subequations}
\begin{align}
	\label{H_subtr}
	(H p_D) (\bd{y}) &=  \int_{\pa\Omega_D} [p_D(\bd{x}) - p_D(\bd{y})] \, \frac{\pa G(\bd{y,x})}{\pa n(\bd{x})} \, dS(\bd{x}) + \frac12 p_D(\bd{y}), \quad \bd{y}\in\pa\Omega_D, \\[6pt]
	\label{K_subtr}
	(K \bd{u}_S) (\bd{y}) &= - \frac{1}{8\pi}\int_{\pa\Omega_S} [\bd{u}_S(\bd{x}) - \bd{u}_S(\bd{y})] \cdot T(\bd{y,x}) \cdot \bd{n}(\bd{x})dS(\bd{x}) - \frac12 \bd{u}_S(\bd{y}), \quad \bd{y}\in\pa\Omega_S.
\end{align}
\end{subequations}

Next, we apply the regularization method of~\cite{beale01, beale04} for Laplace's equation, and recently developed for the Stokes equations in~\cite{tlupova18}. The approach is to replace $1/r^p$, where $r=|\bd{\hat{y}}|$, by a smooth version $s(r/\del)/r^p$, where $s(r/\del)$ is a radial smoothing function and $\del>0$ is a smoothing, or regularization, parameter. As we are interested in evaluating the integrals on the boundary, the smoothing factor $s$ can be chosen to have high order, so that the accuracy in solving the integral equations approaches $O(h^5)$, where $h$ is the grid size in the coordinate planes. The derivation of these smoothing functions can be found in~\cite{beale01, beale04, tlupova18}. For completeness, we state the regularization formulas for both the Darcy and Stokes problems. The kernels in the Darcy representation~\eqref{Greens_bdry} are replaced with their regularized versions,
\beq
	G^\del(\bd{y,x}) = -\frac{s_1(r/\del)}{4\pi r}, \qquad \qquad \nabla G^\del(\bd{y,x}) = \nabla G^\del(\bd{y,x}) s_2(r/\del) = \frac{\bd{\hat{y}}}{4\pi r^3} s_2(r/\del).
\eeq
The regularized single and double layer Stokes kernels in the representation~\eqref{Stokes_bdry} are
\begin{equation}
	S^\del_{ij}(\bd{y,x}) = \del_{ij} \frac{s_1(r/\del)}{r} + \hat{y}_i \ \hat{y}_j \frac{s_3(r/\del)}{r^3}, \qquad \qquad
	T^\del_{ijk} (\bd{y,x}) = -6 \hat{y}_i\hat{y}_j\hat{y}_k \frac{ s_4(r/\del)}{r^5}.
\end{equation}
The radial smoothing functions $s_i(r/\del)$ are chosen to yield high order accuracy when the integrals are evaluated on the surface, with the error approaching $O(h^5)$ in discretization size $h$ as shown in~\cite{beale01, beale04, tlupova18},
\begin{subequations}
\begin{align}
	s_1(r) &= \erf(r) - \frac{2}{3\sqrt{\pi}} r (2 r^2 - 5) e^{-r^2},\\
	s_2(r) &= \erf(r) + \frac{2}{3\sqrt{\pi}} r (2 r^2 - 3) e^{-r^2},\\
	s_3(r) &= \erf(r) - \frac{2}{3\sqrt{\pi}}  r (4 r^4 - 14 r^2 + 3) e^{-r^2}, \\
	s_4(r) &= \erf(r) - \frac{2}{9\sqrt{\pi}}  r (8 r^6 - 36 r^4 + 6 r^2 + 9) e^{-r^2},
\end{align}
\end{subequations}
where $\erf$ is the error function. This approach is quite simple in that it does not require any special quadrature for evaluating the integrals on the boundary.


\subsection{Quadrature}
\label{sec:Quadrature}

\begin{figure}[htb]
\begin{centering}
\scalebox{0.4}{\includegraphics{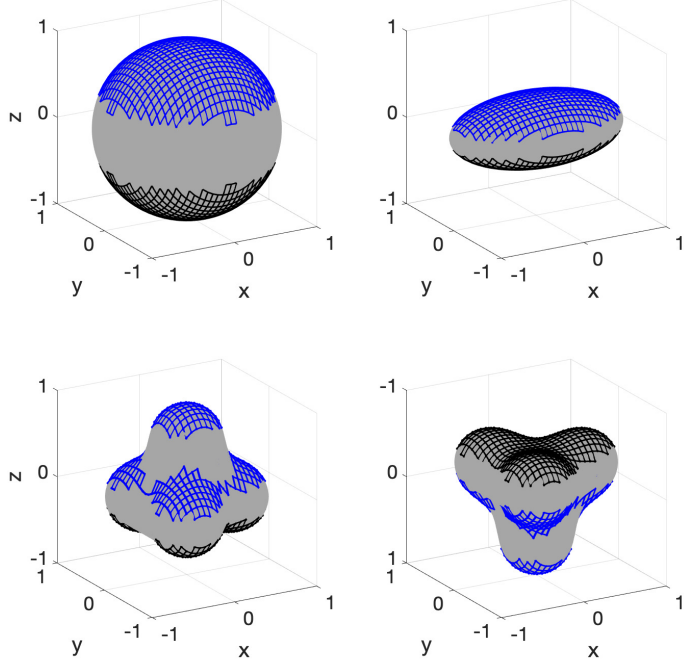}} 
\caption{Quadrature points generated by $\psi_3$ and $h=1/16$ on the surface of the unit sphere~\eqref{Sphere} (top left), the ellipsoid with semiaxes $a=1,b=0.6,c=0.4$~\eqref{Ellipsoid} (top right), and the 'molecule'~\eqref{Molecule} (bottom row). The quadrature points are at the intersections of the lines.}
\label{Surfaces}
\end{centering}
\end{figure}

With the integrands smoothed out, to discretize the integrals a quadrature rule for closed surfaces from~\cite{wilson, beale16} is applied, and we briefly describe it here. First, an angle $\theta$ is chosen and a partition of unity on the unit sphere is defined,
consisting of functions $\psi_1, \psi_2, \psi_3$ with $\Sigma_i \psi_i \equiv 1$ such that
$\psi_i(\bd{n}) = 0$ if $|\bd{n}\cdot\bd{e}_i| \leq \cos{\theta}$, where $\bd{e}_i$ is the
$i$th coordinate vector.  Here we use $\theta = 70^o$.  For mesh size $h$,
a set $R_3$ of quadrature points consists of points $\bf{x}$ on the surface of the form
$(j_1h,j_2h,x_3)$ such that $|\bd{n(x}) \cdot \bd{e}_3| \geq \cos{\theta}$, where
$\bd{n}(\bd{x})$ is the unit normal at $\bd{x}${, see Fig.~\ref{Surfaces}}.  Sets $R_1$ and $R_2$ are
defined similarly.  For a function $f$ on the surface the integral is computed as
\beq
	\label{Quadrature}
	\int_S f(\bd{x}) \,dS(\bd{x}) \approx \sum_{i=1}^3 \sum_{\bd{x} \in R_i} 
            \frac{\psi_i(\bd{n}(\bd{x}))\,f(\bd{x}) }{| \bd{n}(\bd{x})\cdot\bd{e}_i |} \,h^2.
\eeq
{The partition of unity functions $\psi_i$ are constructed from the $C^\infty$ bump function $b(r) = e^{r^2/(r^2 - 1)}$ for $|r| < 1$ and zero otherwise.  The quadrature is effectively reduced to the trapezoidal rule without boundary.  Thus for regular integrands the quadrature has arbitrarily high order accuracy, limited only by the degree of smoothness of the integrand and surface.}  The points in $R_i$
can be found by a line search since they are well separated; see~\cite{wilson} and~\cite{beale16}.

Suppose the surface is given by $\phi(x_1, x_2, x_3) = 0$, with $\phi > 0$ outside, and the normal vector $\nabla\phi$ is predominantly in the $x_3$ direction. Then the parameterization is given by $x_3 = z(x_1, x_2)$, where $z$ is the vertical coordinate on the surface. Derivatives of this parameterization are computed implicitly by $z_i \equiv \pa z/ \pa x_i = -\phi_i/\phi_3$, the tangent vectors are defined as  $\bd{T}_1 = (1, 0, z_1)$ and $\bd{T}_2 = (0, 1, z_2)$. The outward normal is $\bd{n} = \nabla\phi/|\nabla\phi |$ or $\bd{n} = \pm (-z_1,-z_2, 1)/\sqrt{1+z_1^2+z_2^2}$, where "+" is used if $x_3 > z(x_1, x_2)$ outside and "-" is used otherwise.

In Figure~\ref{Surfaces}, we show the set $R_3$ of quadrature points, given by the top and bottom coordinate grids $(j_1h,j_2h,x_3)$, for the following three surfaces:
\begin{subequations}
\begin{align}
	\label{Sphere}
		\phi(x_1,x_2,x_3) &=  x_1^2 + x_2^2 + x_3^2 - 1, \\[6pt]
	\label{Ellipsoid}
		\phi(x_1,x_2,x_3) &= \frac{x_1^2}{a^2} + \frac{x_2^2}{b^2} + \frac{x_3^2}{c^2} - 1, \\[6pt]
	\label{Molecule}
		\phi(x_1,x_2,x_3) &= \sum_{k=1}^4 \exp(- |{\bf x} - {\bf x}_k|^2/r^2) - c.
\end{align}
\end{subequations}

For the ellipsoid~\eqref{Ellipsoid} we set $a=1,b=0.6,c=0.4$, and for the 'four-atom molecule' surface~\eqref{Molecule} we use centers $\bd{x}_1 = (\sqrt{3}/3,0,-\sqrt{6}/12)$, $\bd{x}_{2,3} = (-\sqrt{3}/6,\pm .5,-\sqrt{6}/12)$, $\bd{x}_4 = (0,0,\sqrt{6}/4)$ and $r = .5$, $c = .6$, as in~\cite{beale16}.

This quadrature rule works well for general surfaces without requiring coordinate charts, using projections on coordinate planes instead. As the weights for the quadrature points are found from a partition of unity on the unit sphere, they do not depend on the particular surface. The quadrature points can be found efficiently if, for example, the surface is given analytically or numerically as the level set of a function. With this approach, the formulation is very simple to implement, the data structure and information needed to describe the boundary are minimal, and there are no parameters to fine-tune except the regularization parameter $\delta$ (which is usually taken as a multiple of grid size $h$, see the discussion in~\cite{ tlupova18}).


\section{Algorithm}
\label{sec:Algorithm}

Now we summarize the algorithm for solving the Stokes-Darcy system. First, we note that
\beq
 	\frac{\pa p_D}{\pa n} = -\frac{\mu}{\kappa} \bd{u}_D \cdot \bd{n}_D,
\eeq
and by combining the interface conditions \eqref{BeaversJoseph} and \eqref{Press_contin},
\beq
	\label{f}
	\bd{f}_S = \bdg{\sigma}_S\cdot\bd{n}_S = -p_D \bd{n}_S - \frac{\gamma}{\sqrt{\kappa}} [\bd{u}_S - (\bd{u}_S\cdot \bd{n}_S)\bd{n}_S],
\eeq
so we rewrite the Darcy~\eqref{Darcy_BIEM} and Stokes~\eqref{Stokes_BIEM} integral equations as
\begin{subequations}
\begin{align}
	\label{Darcy_final}
	\frac12 p_D &= H p_D + b_D, \qquad \qquad b_D = \frac{\mu}{\kappa} L (\bd{u}_D \cdot \bd{n}_D),\\[6pt]
	\label{Stokes_final}
	\frac12 \bd{u}_S &= K \bd{u}_S + b_S, \qquad \qquad b_S = -M \bd{f}_S + \bd{u}^\infty,
\end{align}
\end{subequations}

We outline the algorithm of the iterative Stokes and Darcy solutions using successive approximations approach, keeping in mind that the successive approximations are eventually replaced by GMRES for a more efficient solution.
\begin{enumerate}
	\item Let $k=0$, $\bd{u}_S^{(0)} = \bd{0}$ and $q^{(0)} = -\bd{u}_S^{(0)}\cdot \bd{n}_S = 0$.
	\item For $k=1,2,...$ until convergence, do:
	\begin{enumerate}
	\item Compute the Laplace single layer potential $b_D$ in~\eqref{Darcy_final} using $\bd{u}_D\cdot \bd{n}_D = q^{(k-1)}$.
	\item Solve the Darcy problem \eqref{Darcy_final} for $p_D^{(k)}$ via (i) successive approximations \eqref{Darcy_BIEM_iter} or (ii) GMRES.
	\item Compute $\bd{f}_S$ as in~\eqref{f}, using $p_D^{(k)}$ and $\bd{u}^{(k-1)}_S$.
	\item Compute the Stokes single layer potential $b_S$ in~\eqref{Stokes_final}.
	\item Solve the Stokes problem \eqref{Stokes_final} for $\bd{u}_S^{(k)}$ via (i) successive approximations \eqref{Stokes_BIEM_iter} or (ii) GMRES.
	\item Set $q^{(k)} = (1-\theta)q^{(k-1)} - \theta \bd{u}_S^{(k)}\cdot \bd{n}_S$.
	\end{enumerate}
\end{enumerate}
Here $\theta \in (0,1)$ is a relaxation parameter. We assume convergence when the relative error between $q^{(k)}$ and $q^{(k-1)}$, which we call the residual, 
\beq
\label{DDM_residual}
	\left( \frac{\sum_{i=1}^N \left(q^{(k)}_i - q^{(k-1)}_i\right)^2}{\sum_{i=1}^N \left(q^{(k)}_i\right)^2} \right)^{1/2},
\eeq
falls below a prescribed tolerance.

\subsection{Discrete iteration operator}
\label{sec:DiscIterOper}

To analyze the convergence properties of DDM in the discrete solution, we look at
the iteration matrix obtained from the integral formulation and its eigenvalues. For simplicity, we assume $\bd{u}^\infty=0$ and the slip coefficient $\gamma=0$. With the iteration variable defined as $q = -\bd{u}_S\cdot \bd{n}_S$, each Dirichlet-Neumann iteration includes solving the following two problems:
\beq
	\frac12 p_D^{(k)} = \tilde{H}p_D^{(k)} + \frac{\mu}{\kappa} \tilde{L} q^{(k-1)},
\eeq
\beq
	\frac12 \bd{u}_S^{(k)} = \tilde{K} \bd{u}_S^{(k)} + \tilde{M}_n p_D^{(k)},
\eeq
where $\tilde{H}, \tilde{K}$ are the matrices corresponding to the double layer potentials for Laplace and Stokes, respectively, and $\tilde{L}, \tilde{M}_n$ are the matrices corresponding to the single layer for Laplace and Stokes (the Stokes single layer $\tilde{M}_n$ including an extra multiplication by the normal). The 'tilde' symbol denotes the discrete version of the operators, with the kernels smoothed out as in Section~\ref{sec:Regularization} and the quadrature in Section~\ref{sec:Quadrature} applied. All these are dense matrices.

We can then write the discrete representation of the iteration operator as follows:
\begin{eqnarray}
	q^{(k)} &= & (1-\theta)q^{(k-1)} - \theta \bd{u}_S^{(k)}\cdot \bd{n}_S \nonumber \\
	&=& (1-\theta)q^{(k-1)} - \theta \left[(I/2 - \tilde{K} )^{-1} \tilde{M}_n p_D^{(k)} \right]\cdot \bd{n}_S \nonumber \\
	&=& (1-\theta)q^{(k-1)} - \theta \left[(I/2 - \tilde{K} )^{-1} \tilde{M}_n\right] \cdot \bd{n}_S \, p_D^{(k)} \nonumber \\
	&=& (1-\theta)q^{(k-1)} - \theta \left[(I/2 - \tilde{K} )^{-1} \tilde{M}_n\right] \cdot \bd{n}_S \, (I/2-\tilde{H})^{-1} \, \frac{\mu}{\kappa} \, \tilde{L} q^{(k-1)}.
	\label{Discrete_iter_oper}
\end{eqnarray}
We write~\eqref{Discrete_iter_oper} as 
\beq
	q^{(k)} = \tilde{\mathcal{A}} q^{(k-1)},
\eeq
with the iteration matrix 
\beq
\label{A_mat}
	\tilde{\mathcal{A}} = I - \theta \left\{ I + \frac{\mu}{\kappa}\big[(I/2 - \tilde{K} )^{-1} \tilde{M}_n\big] \cdot \bd{n}_S \ \big[(I/2-\tilde{H})^{-1}  \tilde{L}\big] \right\}.
\eeq
The operator in the first square brackets corresponds to solving the Stokes problem for $\bd{u}_S$, while the operator in the second square brackets corresponds to solving the Darcy problem for $p_D$.


\section{Numerical Results}
\label{sec:NumericalResults}

Here we analyze the algorithm in terms of accuracy, convergence, and dependence on physical parameter values for several surfaces. The code was written in C++, compiled with the clang compiler, and all computations were performed on a MacBook Pro 2.3 GHz and 32 GB RAM, running the Big Sur OS. 

\subsection{BIEM accuracy}

We first demonstrate the high achievable accuracy of the integral formulations on the boundary. We solve the Darcy and Stokes integral equations independently of each other, using successive approximations~\eqref{Darcy_BIEM_iter} and~\eqref{Stokes_BIEM_iter}, each with a known exact solution. For the Darcy problem test, we used $\kappa=1$ and the harmonic function $p_D = e^x\sin y$ on the sphere of radius $R=0.8$. For the Stokes problem test, we used the solution similar to Sec. 4.4. of~\cite{tlupova18}. On the outside of a sphere of radius 0.8, we assume the velocity given by a point force singularity of strength $(1,0,0)$, placed at $(0.2, 0,0)$. We then take the solution on the inside of the sphere to be 0, so that the strengths in the Stokes single and double layer potentials are the jumps in velocity and surface force across the surface.
 
\begin{figure}[h]
\begin{centering}
\scalebox{0.4}{\includegraphics{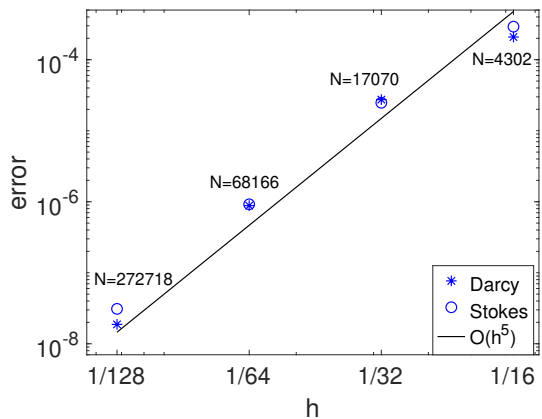}} 
\caption{Convergence in grid size $h$ of separate Darcy and Stokes problems solved on the surface of a unit sphere by successive approximations \eqref{Darcy_BIEM_iter} and \eqref{Stokes_BIEM_iter}, with $\delta=3h$.}
\label{Errors}
\end{centering}
\end{figure}
Figure~\ref{Errors} shows the accuracy of this solution, where error is defined as the $L_2$ norm of the difference in the computed and exact solutions, $p_D$ for Darcy and $\bd{u}_S$ for Stokes. We used a smoothing parameter $\delta = 3h$ here and throughout this paper, where $h$ is the grid size in the coordinate planes. A larger choice of $\delta$ is needed to ensure the regularization error is dominant over the discretization error, so that the total error approaches $O(h^5)$ (see~\cite{tlupova18}). As expected, the convergence rate observed in Fig.~\ref{Errors} is $O(h^5)$. Other surfaces, such as a thin ellipsoid, are generally expected to give somewhat larger errors due to the larger curvature and varied spacing~\cite{tlupova18}.


\subsection{Discrete iteration operator}

\begin{figure}[htb]
\begin{centering}
	\scalebox{0.27}{\includegraphics{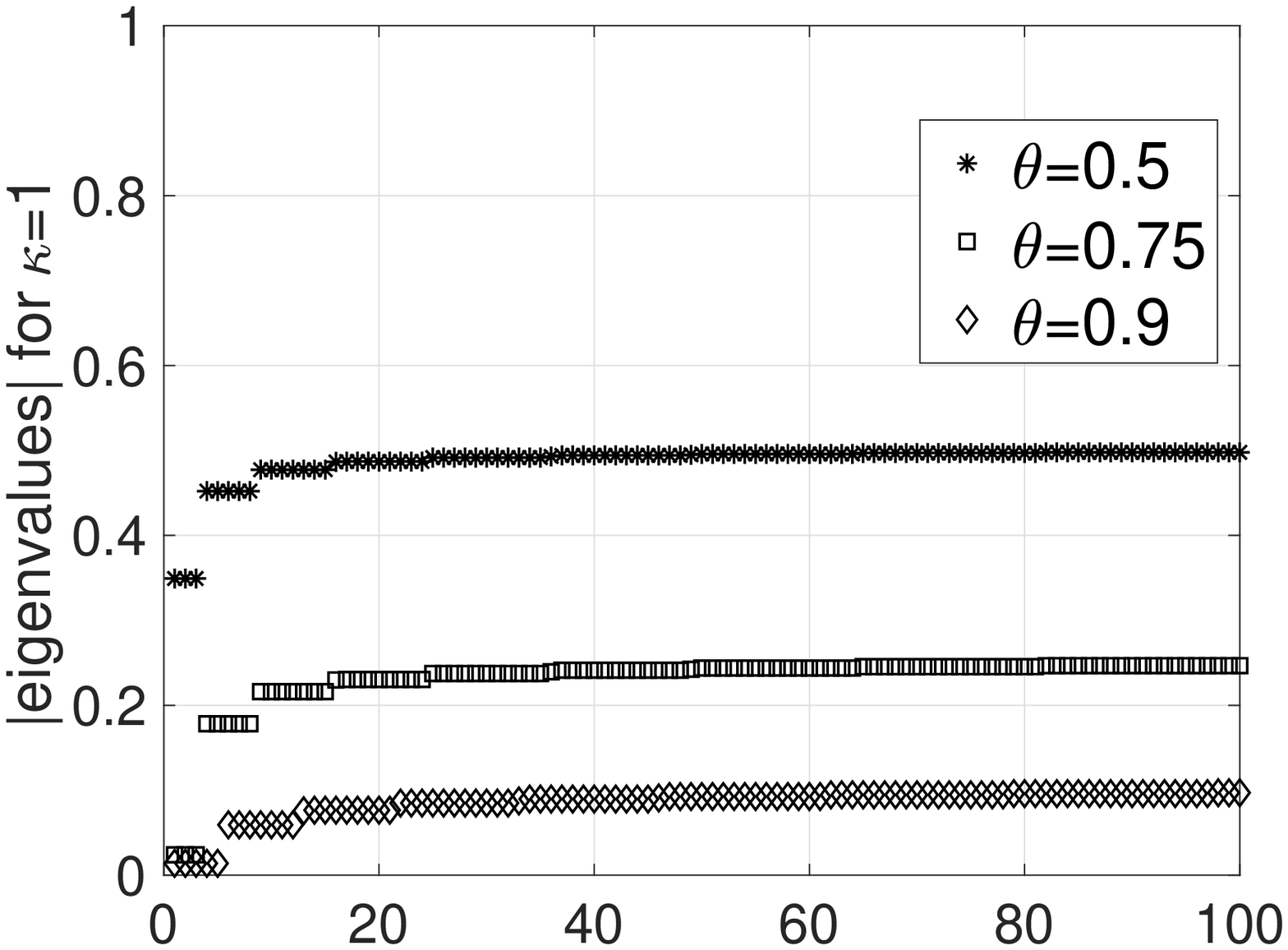}}
	\scalebox{0.27}{\includegraphics{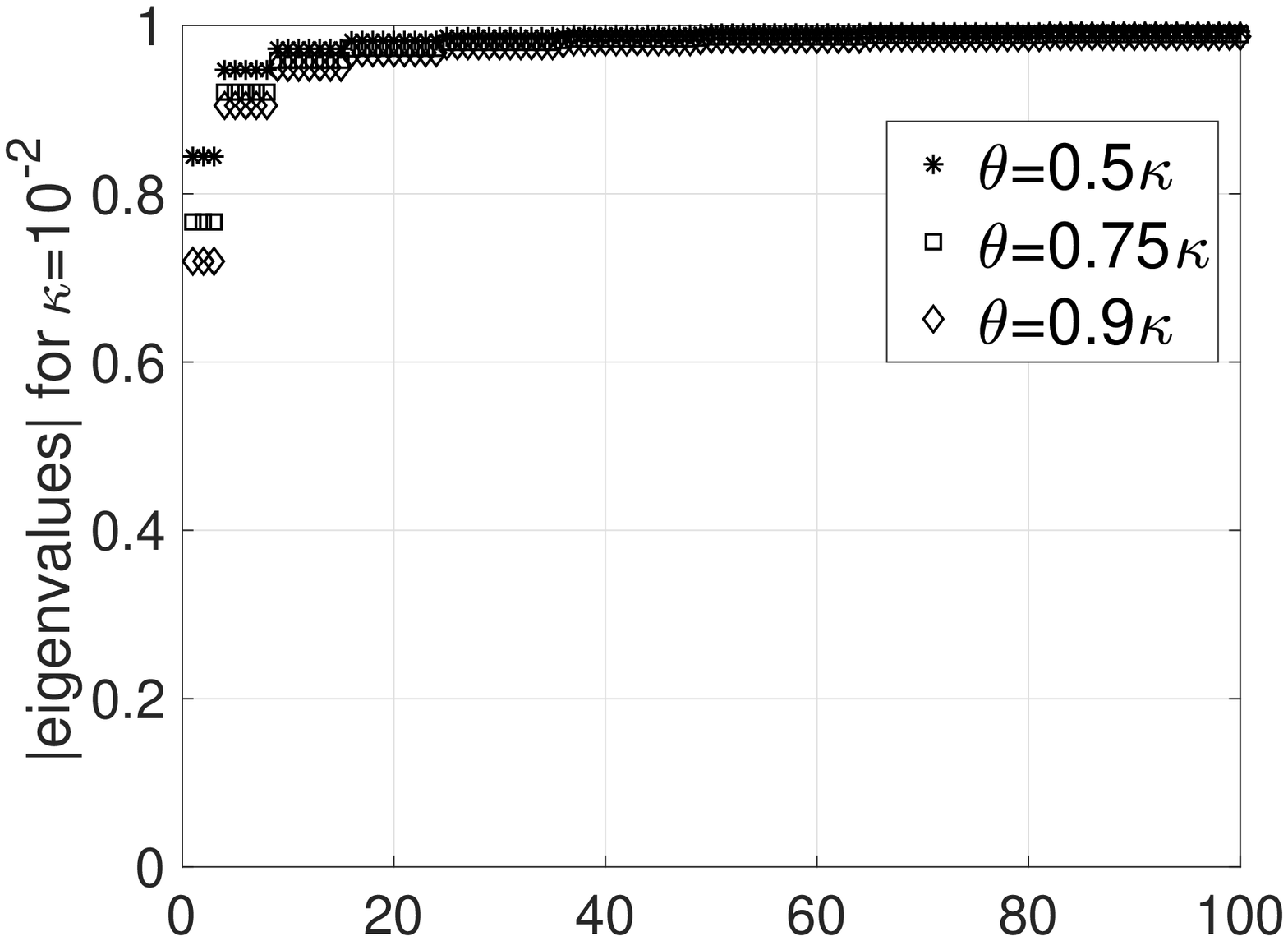}}
	\scalebox{0.27}{\includegraphics{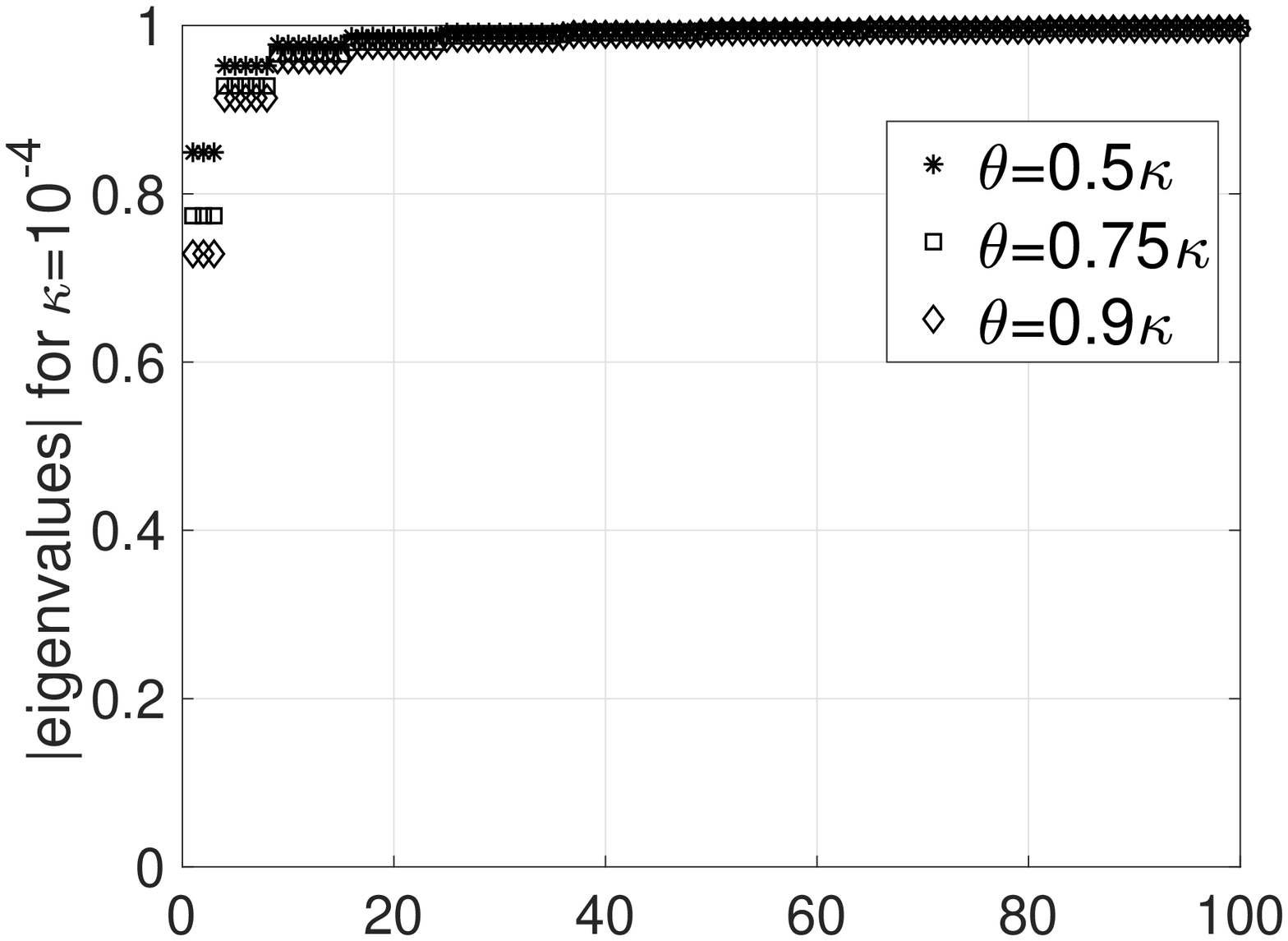}}
	\caption{Eigenvalues of the iteration matrix $\tilde{\mathcal{A}}$ in~\eqref{A_mat} on the unit sphere for $h=1/16$, for permeability values $\kappa=1, 10^{-2}, 10^{-4}$.}
\label{Sphere_eigenvalues}
\end{centering}
\end{figure}

Here we validate the numerical implementation of the Dirichlet-Neumann iterations for the spherical geometry. Similar to the setup used in Section~\ref{sec:Analysis}, consider the Darcy problem inside the unit sphere centered at the origin, coupled with the Stokes problem outside the sphere. Parameters are set to $\mu=1$, $\gamma=0$, and the unit sphere was discretized using $h=1/16$ leading to $N=4302$ quadrature points. In Figure~\ref{Sphere_eigenvalues}, the eigenvalues of the discrete iteration matrix $\tilde{\mathcal{A}}$ in~\eqref{A_mat} are computed and the first 50 are shown, to compare with the harmonic coefficients $\mathcal{A}_n$ of the iteration operator in Figure~\ref{Sphere_modes}. 
\begin{itemize}
	\item $\kappa=1$. In the left figure, the eigenvalues are shown for permeability $\kappa=1$ and varying relaxation parameter $\theta$ values. Again we see that a larger value of $\theta \in (0,1)$ leads to smaller eigenvalues, where most cluster around the value $1-\theta$. We then expect the successive approximations to converge faster for larger $\theta$. 
	\item $\kappa \ll 1$. Next, keeping $\mu=1$ fixed, we compute the eigenvalues for smaller values of permeability, $\kappa=10^{-2}$ (center figure), and $\kappa=10^{-4}$ (right figure). For these, we set the relaxation parameter $\theta=c\kappa$, and vary the value of coefficient $c=0.5, 0.75, 0.9$ as before. In all cases, the vast majority of the eigenvalues cluster around 1 (we saw this with the harmonic coefficients as well). As we will show, this is going to severely inhibit convergence speeds of DDM by successive approximations, and we will perform the DDM iterations using GMRES instead.
\end{itemize}
	
\subsection{Benchmark problem: Stokes flow past a porous sphere}

In this section we test our algorithm using the problem of uniform flow of a viscous fluid past a porous spherical shell. Assuming that the fluid in the permeable sphere $0\leq r \leq R$ obeys Darcy's law and the fluid outside the sphere satisfies the Stokes equations, the solution was obtained by Joseph and Tao~\cite{joseph-tao64}. The matching conditions enforced on the sphere are continuity of pressure and normal velocity, as well as the no-slip condition for the tangential velocity of the exterior flow. Considering a stationary porous sphere of radius $R$ in a uniform streaming viscous fluid $\bd{U} = (0,0,U)$, the solution is given in spherical coordinates $(r,\theta,\phi)$,
\begin{eqnarray}
	p^D &= &-3\mu\, U\, r\, \cos\theta / (2R^2+\kappa),\\
	u^D_r &=& 3 U\, \kappa\, \cos\theta / (2R^2+\kappa),\\
	u^D_\theta &=& -3 U\, \kappa\, \sin\theta / (2R^2+\kappa),
\end{eqnarray}
for the Darcy problem, and 
\begin{eqnarray}
	p^S &= &-3\mu\, R\, U\, \cos\theta / r^2(2+\kappa/R^2),\\
	u^S_r &=& \left\{ \frac{-3 R U}{(2+\kappa/R^2)\, r} \left( 1-\frac{R^2}{3r^2} \left(1+\frac{2\kappa}{R^2} \right) \right)  + U \right\} \cos\theta,\\
	u^S_\theta &=& \left\{ \frac{3 R U}{(2+\kappa/R^2)\, 2\, r} \left( 1+\frac{R^2}{3r^2} \left(1+\frac{2\kappa}{R^2} \right) \right)  - U \right\} \sin\theta,
\end{eqnarray}
for the Stokes problem, where $u_\phi=0$ due to axial symmetry. The drag on the porous sphere, defined as the hydrodynamic force exerted on the sphere,
\beq
\label{drag_int}
	\bd{D} = \int_{\pa\Omega} \bd{f}\, dS(\bd{x}),
\eeq
where $\bd{f} = \bdg{\sigma}_S\cdot\bd{n}$, $\bd{n}=\bd{e}_r$, was found to equal
\beq
\label{drag_value}
	\bd{D} = \frac{6\pi\, \mu\, R\, U}{1+\kappa/(2R^2)}.
\eeq
The usual solution and drag for the flow past a solid sphere are recovered from these formulas by setting $\kappa=0$.

To test our algorithm against this solution, we compute the Stokes stress vector as $\bd{f} = f_r \bd{e}_r + f_\theta \bd{e}_\theta$ with
\beq
\label{f_r}
	f_r = -p_D + 2\mu \frac{\pa u_r}{\pa r},
\eeq
\beq
\label{f_theta}
	f_\theta = \mu \left( r \frac{\pa}{\pa r} \left(\frac{u_\theta}{r}\right) + \frac{1}{r} \frac{\pa u_r}{\pa\theta} \right),
\eeq
being the normal and tangential components. Note that we replaced $p_S$ with $p_D$ in~\eqref{f_r} since the exact solution was derived assuming continuous pressure $p_D=p_S$ on the boundary. To compute the drag numerically, we integrate~\eqref{drag_int} using our usual quadrature, and compare to the exact value given by~\eqref{drag_value}.

\begin{figure}[!htb]
\begin{centering}
	\scalebox{0.4}{\includegraphics{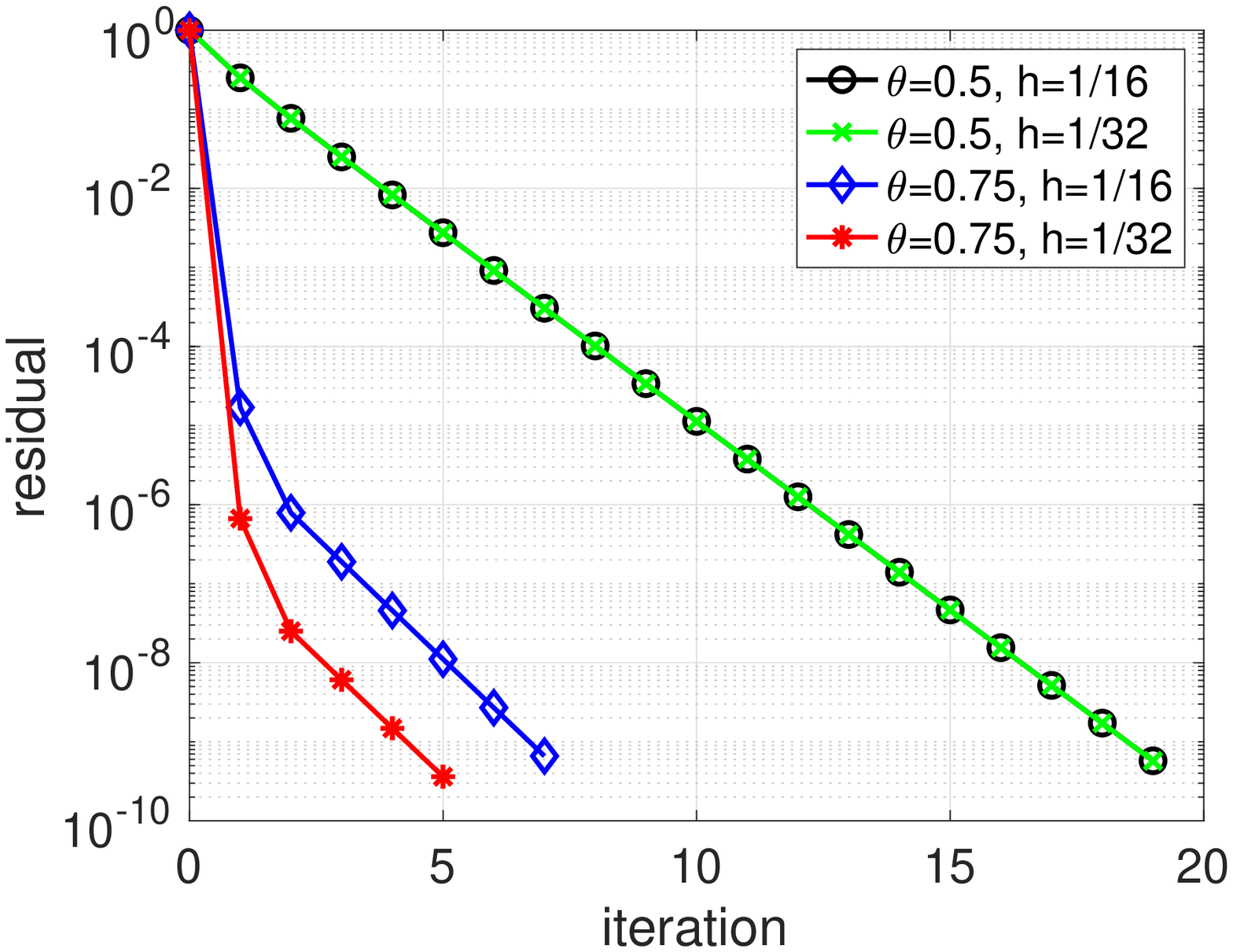}}
	\scalebox{0.4}{\includegraphics{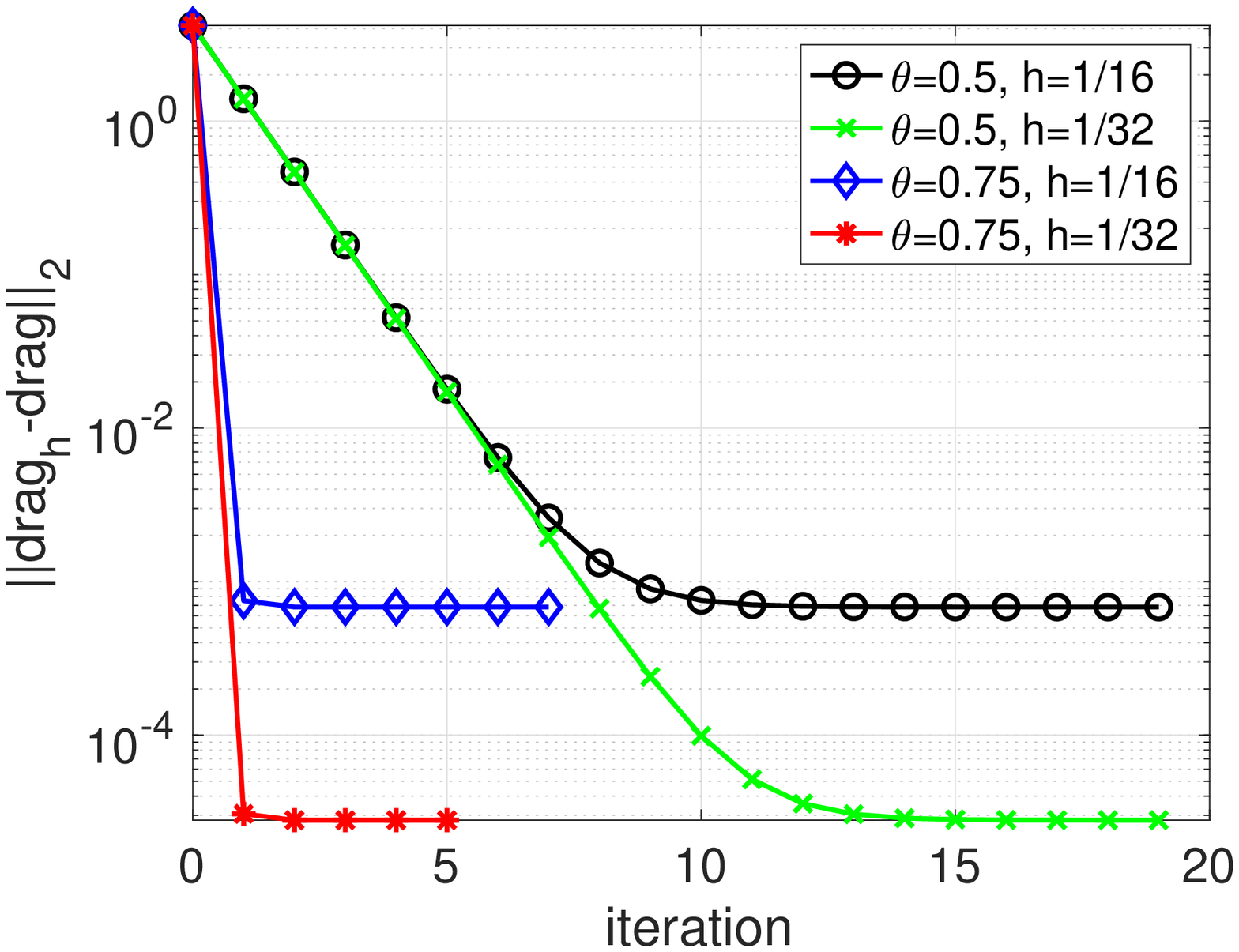}} \\
	\scalebox{0.4}{\includegraphics{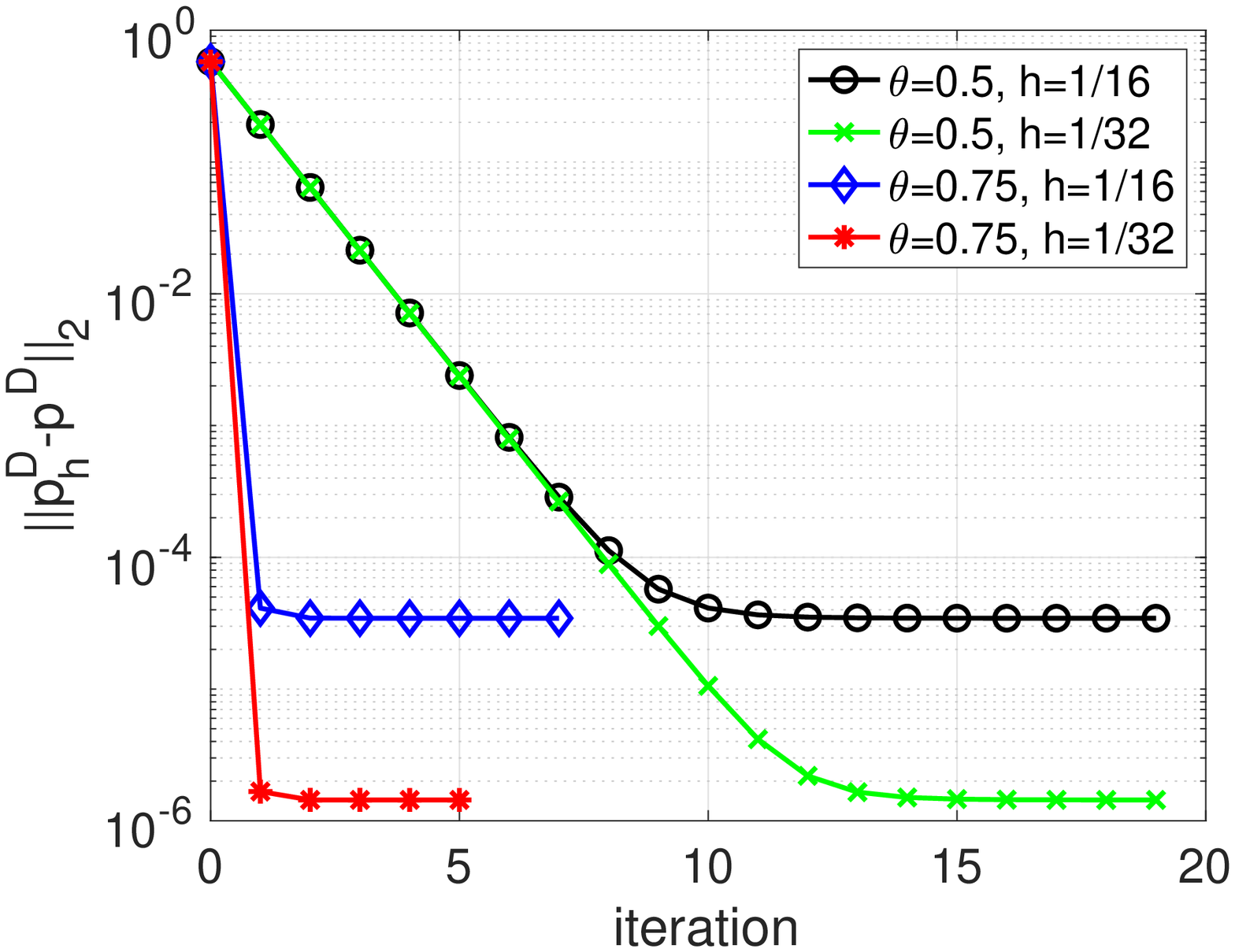}}
	\scalebox{0.4}{\includegraphics{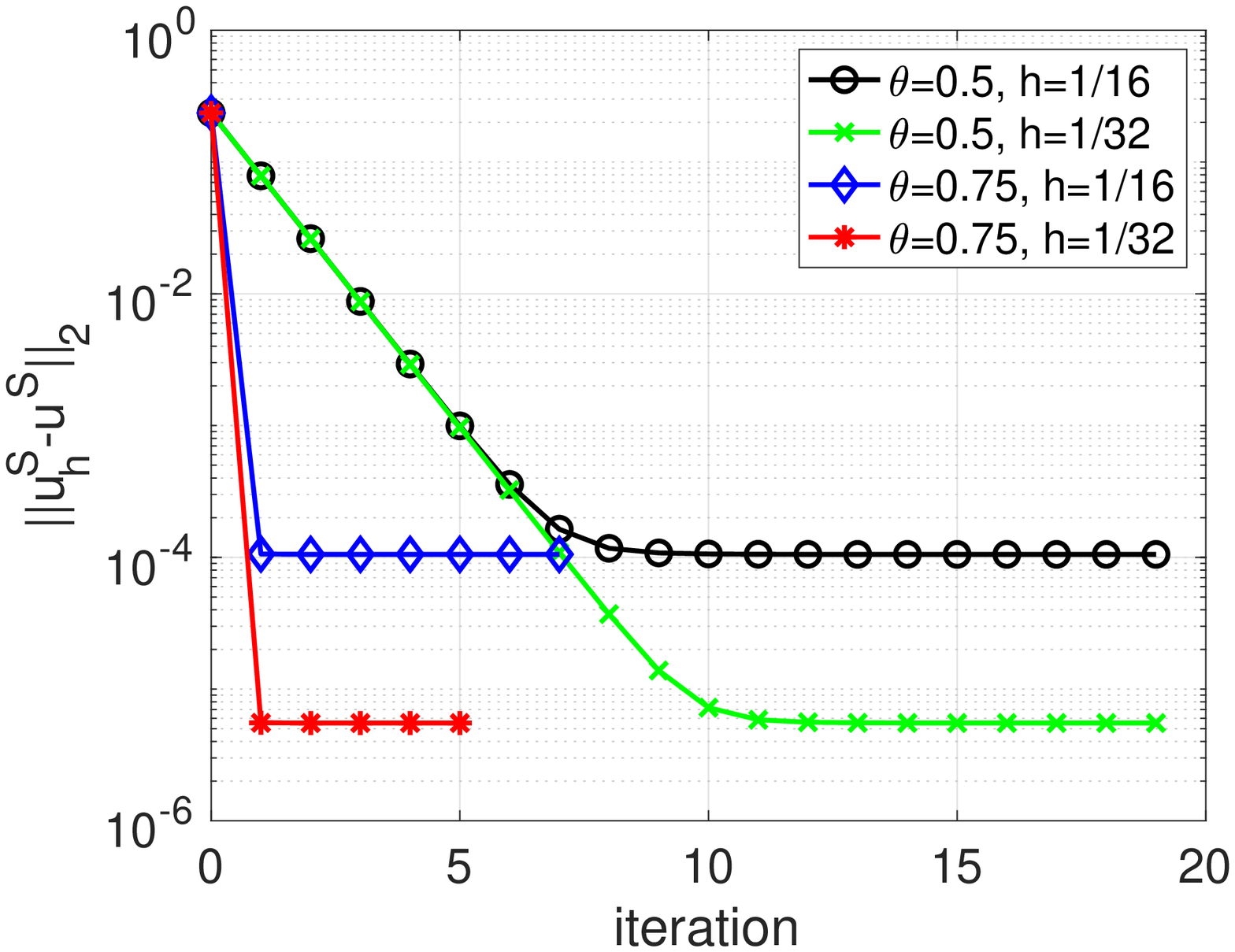}} 
	\caption{Flow past a porous sphere with $\kappa=1$, tolerance $10^{-9}$ for residual/local solutions. Top left: DDM residual, top right: error in drag, bottom left: error in $p_D$, bottom right: error in $\bd{u}_S$.}
\label{Sphere_k1_resi_err}
\end{centering}
\end{figure}

In our tests, we use a unit sphere $R=1$, $U=1, \mu=1$ and no-slip $\gamma=0$. We use a tolerance of $10^{-9}$ for the residual~\eqref{DDM_residual} in the Dirichlet-Neumann iterations, and for the solution of the integral equations, we use successive approximations~\eqref{Darcy_BIEM_iter} and~\eqref{Stokes_BIEM_iter}, as well as GMRES, both with tolerance $10^{-9}$. Table~\ref{table:Sphere_k1} shows the results for $\kappa=1$. For different grid sizes $h$, displayed are: the number of quadrature points $N$ on the sphere surface, the number of Dirichlet-Neumann iterations via successive approximations (D-N SA) for two values of $\theta=0.5$ and $\theta=0.75$, the number of successive approximations (local SA) needed to solve the Darcy and Stokes problems at each D-N step to reach the set tolerance, the number of GMRES iterations (local GMRES) for each problem at each D-N step for the same tolerance, and finally the errors in the Darcy and Stokes solutions computed in $L_2$ norm as
\beq
\label{error_Darcy}
	||p^D_h-p^D||_2 = \left(\sum_{i=1}^N \left(p_D(\bd{x}_i)- p_D^{ex}(\bd{x}_i)\right)^2 \Big/ N \right)^{1/2},
\eeq
\beq
\label{error_Stokes}
	||\bd{u}^S_h-\bd{u}^S||_2 = \left(\sum_{i=1}^N \left|\bd{u}_S(\bd{x}_i)- \bd{u}_S^{ex}(\bd{x}_i)\right|^2 \Big/ N \right)^{1/2},
\eeq
where $p_D^{ex}$ and $\bd{u}_S^{ex}$ are the exact values and $|\bd{u}|$ is the Euclidean norm. For a more detailed view, Figure~\ref{Sphere_k1_resi_err} shows the DDM residual and the errors in drag, $p_D$ and $\bd{u}_S$, by comparing the numerical solution to the exact solution at every D-N iteration for grid sizes $h=1/16$ and $h=1/32$ and varying $\theta$. 

\begin{figure}[!htb]
\begin{centering}
	\scalebox{0.4}{\includegraphics{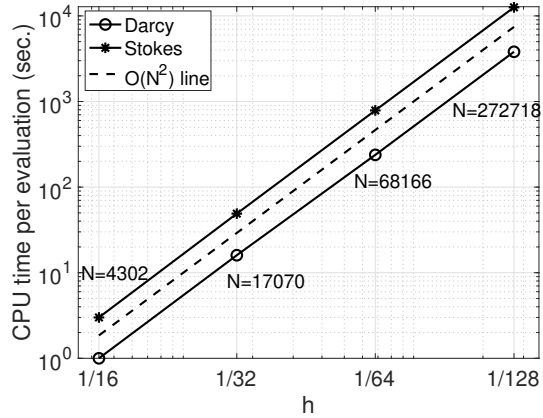}} 
	\caption{Flow past a porous sphere with $\kappa=1$. CPU time in seconds on each evaluation of the double layer potential (one successive approximation or one GMRES iteration).}
\label{Sphere_k1_CPU}
\end{centering}
\end{figure}

We make the following observations:
\begin{itemize}
	\item \underline{D-N convergence and dependence on the relaxation parameter $\theta$}. For $\theta=0.5$, about 19 D-N iterations are required to reach the tolerance, independent of the grid size, while for $\theta=0.75$ only 7 iterations are needed for $h=1/16$, and this number reduces quickly to 2 iterations for smaller grid sizes. Figure~\ref{Sphere_k1_resi_err} (top left) shows the gradual decrease in the residual for $\theta=0.5$ and a more rapid decrease for $\theta=0.75$. This faster convergence can be attributed to the smaller spectral radius of the iteration operator, as was observed for $\kappa=1$ in Figures~\ref{Sphere_modes} and~\ref{Sphere_eigenvalues}. This is remarkable given that with smaller $h$ the local Stokes and Darcy problems become larger.
	\item \underline{Solution errors}. Further, we see that the errors in drag, $p_D$ and $\bd{u}_S$ in Figure~\ref{Sphere_k1_resi_err} all level off at the BIEM accuracy shown in Table~\ref{table:Sphere_k1} (again, much faster with $\theta=0.75$). From the error values in Table~\ref{table:Sphere_k1} we see that the convergence rate is $O(h^5)$ or higher with grid refinement. 
	\item \underline{Local problems: SA vs. GMRES}. Table~\ref{table:Sphere_k1} shows the number of successive approximations (local SA) performed at each D-N iteration, as well as the number of GMRES iterations (local GMRES) as an alternative. For SA, a range is shown because the number decreases slowly with each D-N iteration. It is clear that successive approximations do not converge quickly, and GMRES is far superior in solving the local Stokes and Darcy problems at each iteration. Both methods require evaluating the double layer potential at each iteration, and since GMRES performs far fewer double layer evaluations, we use GMRES in all our remaining tests. 
	\item \underline{Local problems: CPU time}. Both SA and GMRES require evaluating the double layer potential at each iteration, which in discretized form is a dense linear system. The Darcy system is $N\times N$, while the Stokes system is $3N\times 3N$, where $N$ is the number of quadrature points. The computational cost is then $O(N^2)$ for each problem. The actual CPU times for each double layer integral evaluation are shown in Figure~\ref{Sphere_k1_CPU}. Note that the discretized integrals were implemented to take advantage of the target-source symmetry in all kernels, reducing the computational time in half. If very large systems are considered, the computational time could be reduced to $O(N\log N)$ or $O(N)$ using a fast summation algorithm such as a treecode~\cite{wang-krasny-tlupova2} or a fast multipole method~\cite{tornberg-greengard-08}.
	\item \underline{Dependence on $\kappa$}. We next proceed to study the convergence of D-N iterations for different values of permeability. For $\kappa=1$, we see a strong dependence on the values of $\theta$ when successive approximations are used, because the spectral radius of the iteration operator is smaller for larger $\theta$. Table~\ref{table:Sphere_diff_k} reiterates some of these results and also shows the results for $\kappa=10^{-2}$ and $\kappa=10^{-4}$. For small values of $\kappa$ such as these, we saw that the spectral radius of the iteration operator is very close to 1, and this indeed inhibits convergence. The number of successive approximations needed for D-N convergence (D-N by SA in Table~\ref{table:Sphere_diff_k}) exceeds 100 because the residual decreases very slowly. This is not improved by using a larger $\theta$ as it did for $\kappa=1$ (in all cases, we set $\theta=c\kappa$ with $c\in (0,1)$).  Because of this, we perform the D-N iterations using GMRES, which converges much faster. The D-N by GMRES column in Table~\ref{table:Sphere_diff_k} shows the number of iterations for the same tolerance $10^{-9}$ for the residual. The number of iterations is 7 for $\kappa=10^{-2}$ and 8 for $\kappa=10^{-4}$ for the coarsest grid $h=1/16$. As the grid size is reduced and the solution becomes more accurate, the number of D-N iterations by GMRES reduces. For the solution errors, which are also shown in the left graph in Figure~\ref{Errors_All}, we observe convergence of the order $O(h^5)$.
\end{itemize}

\begin{table}[!htb]
  \begin{center}
    \begin{tabular}{|c|c||c|c||c|c||c|c||c|c|}
     \multirow{2}{*}{$h$} & \multirow{2}{*}{$N$} & \multicolumn{2}{c||}{D-N SA} & \multicolumn{2}{c||}{local SA} & \multicolumn{2}{c||}{local GMRES} & \multirow{2}{*}{$||p^D_h-p^D||_2$} & \multirow{2}{*}{$||\bd{u}^S_h-\bd{u}^S||_2$}\\
     & & $\theta=0.5$ & $\theta=0.75$ & Darcy & Stokes & Darcy & Stokes & & \\
      \hline
      1/16 & 4302 & 19 & 7 & 37-2 & 26-2 & 5 & 8 & 3.450e-05 & 1.053e-04 \\
      1/32 & 17070 & 19 & 5 & 37-1 & 26-1 & 4 & 7 & 1.442e-06 & 5.525e-06 \\
      1/64 & 68166 & - & 2 & - & - & 3 & 5 & 2.605e-08 & 1.716e-07 \\
      1/128 & 272718 & - & 1 & - & - & 2 & 3 & 3.236e-10 & 4.735e-09 \\
    \end{tabular}
    \caption{Flow past a porous sphere with $\kappa=1$, DDM residual tolerance = $10^{-9}$, local solution tolerance = $10^{-9}$.}
    \label{table:Sphere_k1}
  \end{center}
\end{table}

\begin{table}[!htb]
  \begin{center}
    \begin{tabular}{|c|c||c|c||c|c||c|c||c|c|}
     \multirow{2}{*}{$\kappa$} & \multirow{2}{*}{$h$} & \multicolumn{2}{c||}{D-N by SA} & \multicolumn{2}{c||}{D-N by GMRES} & \multicolumn{2}{c||}{local GMRES} & \multirow{2}{*}{$||p^D_h-p^D||_2$} & \multirow{2}{*}{$||\bd{u}^S_h-\bd{u}^S||_2$} \\
     & & $\theta=0.5\kappa$ & $\theta=0.75\kappa$ & $\theta=0.5\kappa$ & $\theta=0.75\kappa$ & Darcy & Stokes & & \\
      \hline
      $1$ & 1/16 & 19 & 7 & 4 & 4 & 5 & 8 & 3.450e-05 & 1.053e-04 \\
             & 1/32 & 19 & 5 & 4 & 4 & 6 & 6 & 1.442e-06 & 5.525e-06 \\
             &1/64 & - & 2 & 3 & 3 & 6 & 6 & 2.605e-08 & 1.716e-07 \\
      \hline
      $10^{-2}$ & 1/16 & $>100$ & $>100$ & 7 & 7 & 6 & 9 & 9.484e-05  & 4.500e-05 \\
                      & 1/32 & $>100$ & $>100$ & 6 & 7 & 7 & 7 & 3.498e-06 & 2.576e-06 \\
                      &1/64 & - & - & 4 & 4 & 8 & 7 & 3.084e-08 & 8.546e-08 \\
      \hline
      $10^{-4}$ & 1/16 & $>100$ & $>100$ & 8 & 8 & 8 & 11 & 5.474e-04 & 4.243e-05 \\
                      & 1/32 & $>100$ & $>100$ & 2 & 2 & 8 & 9 & 3.490e-06 & 2.554e-06 \\
                      &1/64 & - & - & 1 & 1 & 5 & 5 & 3.899e-08 & 8.451e-08 \\
    \end{tabular}
    \caption{Flow past a porous sphere with varying $\kappa$, tolerance $10^{-9}$ for residual/local solutions.}
    \label{table:Sphere_diff_k}
  \end{center}
\end{table}


\subsection{Thin porous ellipsoid}

Next we test the method on a geometry that has a larger curvature, a thin ellipsoid defined in~\eqref{Ellipsoid}
with $a=1,b=0.6,c=0.4$. In this case the exact solution is not known, and we compute the error empirically as $||p_h - p_{h/2}||_2$ and $||\bd{u}_h - \bd{u}_{h/2}||_2$, with similar definitions of the $L_2$ norm as in~\eqref{error_Darcy}-\eqref{error_Stokes}. Table~\ref{table:Ellipsoid_diff_k} shows the results for permeability values $\kappa = 1, 10^{-2}, 10^{-4},10^{-7}$. As we saw in the sphere example, the number of DDM iterations increases slightly when we go from $\kappa=1$ to $\kappa<1$. Unlike the sphere, however, as the grid is refined, the number of DDM iterations does not decrease but rather increases slightly. Furthermore, for the small value of the permeability $\kappa=10^{-4}$, the DDM iterations do not converge to the prescribed tolerance $10^{-9}$ quickly. In the number of iterations shown in the table (denoted by *), the relative residual reaches a tolerance of about $10^{-6}$, with the consecutive residuals differing by less than $10^{-9}$. As expected for the thin ellipsoid, the solution accuracy is lower than what we observed with the sphere.

\begin{table}[!htb]
  \begin{center}
    \begin{tabular}{|c|c|c||c|c||c|c||c|c|}
     \multirow{2}{*}{$\kappa$} & \multirow{2}{*}{$h$} & \multirow{2}{*}{$N$} & \multicolumn{2}{c||}{D-N by GMRES} & \multicolumn{2}{c||}{local GMRES} & \multirow{2}{*}{$||p^D_h-p^D_{h/2}||_2$} & \multirow{2}{*}{$||\bd{u}^S_h-\bd{u}^S_{h/2}||_2$} \\
     & & & $\theta=0.5\kappa$ & $\theta=0.75\kappa$ & Darcy & Stokes & & \\
      \hline
      $1$ & 1/16 & 1742 & 5 & 5 & 7 & 11 & 3.731e-04 & 6.052e-04 \\
             & 1/32 & 6902 & 5 & 5 & 6 & 10 & 3.486e-05 & 4.399e-05 \\
             &1/64 & 27566 & 5 & 5 & 5 & 9 & 3.300e-06 & 3.589e-06 \\
             &1/128 & 110250 & 5 & 5 & 5 & 8 & - & - \\
      \hline
      $10^{-2}$ & 1/16 & 1742 & 9 & 9 & 9 & 14 & 2.378e-02 & 5.744e-03 \\
                      & 1/32 &  6902 & 12 & 12 & 9 & 14 & 1.772e-03 & 1.797e-03 \\
                      &1/64 & 27566 & 12 & 12 & 8 & 11 & 8.347e-05 & 8.423e-05 \\
                      &1/128 & 110250 & 12 & 12 & 8 & 11 & - & - \\
      \hline
      $10^{-4}$ & 1/16 & 1742 & 9* & 8* & 11 & 16 & 1.689e-01 & 1.873e-02 \\
                      & 1/32 & 6902 & 11* & 11* & 11 & 17 & 2.913e-02 & 2.914e-02 \\
                      &1/64 & 27566 & 13* & 13* & 11 & 17 & 1.435e-03 & 1.436e-03 \\
                      &1/128 & 110250 & 13* & 13* & 11 & 18 & - & - \\
      \hline
      $10^{-7}$ & 1/16 &  1742 & 
       6 &  6 &  14 &  20 & 
       1.900e-01 &  1.865e-02 \\
                      &  1/32 &  6902 & 
                       8 &  8 &  11 &  20 & 
                       3.406e-02 &  3.408e-02 \\
                      & 1/64 &  27566 & 
                       10 &  10 &  11 &  18 & 
                       1.448e-03 &  1.449e-03 \\
                      & 1/128 &  110250 & 
                       10 &  10 &  12 &  18 & - & - \\
    \end{tabular}
    \caption{Flow past a porous ellipsoid with varying $\kappa$, tolerance $10^{-9}$ for residual/local solutions.}
    \label{table:Ellipsoid_diff_k}
  \end{center}
\end{table}


\subsection{Porous molecule}

\begin{figure}[!htb]
	\hspace{-1.2cm} \scalebox{0.4}{\includegraphics{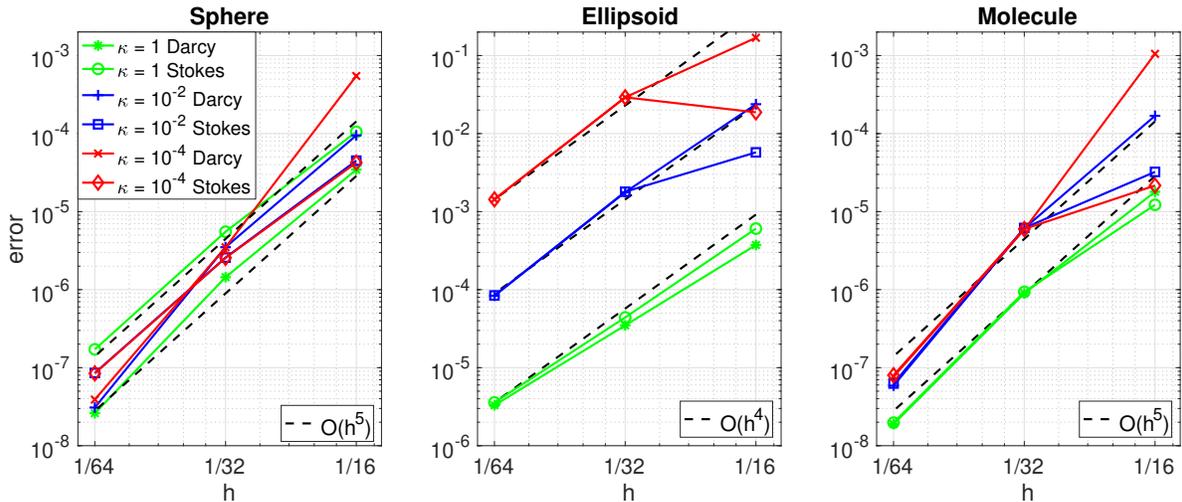}} 
	\caption{Solution errors for three surfaces: unit sphere, thin ellipsoid, four-atom molecule.}
\label{Errors_All}
\end{figure}

The four-atom molecule surface used in these tests is defined in~\eqref{Molecule}. Again, we compute the errors empirically similar to the ellipsoid. The results for different permeability values and grid sizes are shown in Table~\ref{table:Molecule_diff_k}. The method performs well for this surface, converging to the tolerance in just a few iterations, with the speed of convergence improving when the discretization of the surface is refined. For $h=1/128$ for example, which results in $N=272718$ quadrature points, only 2 DDM iterations are needed for $\kappa=1$ and $\kappa=10^{-2}$, while the $\kappa=10^{-4}$ { and $\kappa=10^{-7}$ cases only require} one iteration for the tolerance of $10^{-9}$. The solution accuracy is higher than for the ellipsoid and similar to the sphere accuracy, with the order of convergence as predicted $O(h^5)$ - this is shown in Figure~\ref{Errors_All}.

\begin{table}[!htb]
  \begin{center}
    \begin{tabular}{|c|c|c||c|c||c|c||c|c|}
     \multirow{2}{*}{$\kappa$} & \multirow{2}{*}{$h$} & \multirow{2}{*}{$N$} & \multicolumn{2}{c||}{D-N by GMRES} & \multicolumn{2}{c||}{local GMRES} & \multirow{2}{*}{$||p^D_h-p^D_{h/2}||_2$} & \multirow{2}{*}{$||\bd{u}^S_h-\bd{u}^S_{h/2}||_2$} \\
     & & & $\theta=0.5\kappa$ & $\theta=0.75\kappa$ & Darcy & Stokes & & \\
      \hline
      $1$ & 1/16 & 4302 & 4 & 4 & 5 & 7 & 1.793e-05 & 1.221e-05 \\
             & 1/32 & 17070 & 4 & 4 & 6 & 6 & 8.947e-07 & 9.358e-07 \\
             &1/64 & 68166 & 3 & 3 & 6 & 6 & 1.892e-08 & 1.974e-08 \\
             &1/128 & 272718 & 2 & 2 & 6 & 6 & - & - \\
      \hline
      $10^{-2}$ & 1/16 & 4302 & 8 & 8 & 6 & 9 & 1.689e-04 & 3.228e-05 \\
                      & 1/32 & 17070 & 7 & 7 & 7 & 7 & 6.102e-06 & 6.174e-06 \\
                      &1/64 & 68166 & 4 & 5 & 7 & 7 & 5.877e-08 & 6.263e-08 \\
                      &1/128 & 272718 & 2 & 2 & 7 & 7 & - & - \\
      \hline
      $10^{-4}$ & 1/16 & 4302 & 8 & 11 & 8 & 10 & 1.052e-03 & 2.158e-05 \\
                      & 1/32 & 17070 & 2 & 8 & 9 & 9 & 5.883e-06 & 5.981e-06 \\
                      &1/64 & 68166 & 1 & 1 & 5 & 5 & 7.676e-08 & 8.000e-08 \\
                      &1/128 & 272718 & 1 & 1 & 3 & 3 & - & - \\
      \hline
       $10^{-7}$ &  1/16 &  4302 & 
       1 &  1 &  8 &  12 &
        1.330e-04 &  4.356e-05 \\
                      &  1/32 &  17070 & 
                       1 &  1 &  8 &  10 & 
                       5.710e-06 &  5.797e-06 \\
                      & 1/64 &  68166 & 
                       1 &  1 &  7 &  8 & 
                       7.679e-08 &  8.002e-08 \\
                      & 1/128 &  272718 &  1 &  1 &  3 &  3 & - & - \\
    \end{tabular}
    \caption{Flow past a porous molecule with varying $\kappa$, tolerance $10^{-9}$ for residual/local solutions.}
    \label{table:Molecule_diff_k}
  \end{center}
\end{table}

\section{Conclusions and future work}

We have presented a domain decomposition method of a sequential Dirichlet-Neumann type to solve the Stokes-Darcy system of PDEs. The method is based on splitting the interface conditions so that first, the Darcy problem is solved with one of the conditions, and then the Stokes problem is solved with the other two. The information is exchanged at each step of this iterative procedure. The local Stokes and Darcy problems are solved by second kind boundary integral formulations. A regularization technique specifically designed to achieve high accuracy on the boundary is applied, and the quadrature method of \cite{wilson, beale16} is used, which does not use coordinate charts and is simple to implement. Numerical results validate the algorithm against a problem with a known analytical solution of viscous flow past a porous sphere, as well as demonstrate the applicability of the method to more general geometries. 

Our convergence analysis of the method shows strong dependence of the spectral radius of the iteration operator on the permeability $\kappa$. For small values of $\kappa$, the eigenvalues cluster close to 1, inhibiting convergence of successive approximations. While we demonstrate that using a Krylov subspace method such as GMRES mitigates this issue, more robust DDMs that are less dependent on the physical parameters are desirable and will be investigated. It was shown in the two-dimensional case~\cite{boubendir-tlupova-13} that a Robin-Robin type domain decomposition method with non-local operators (based on the boundary integral formulations) can be chosen in the de-coupling of transmission conditions to dramatically improve the convergence properties of the solution. The dependence of convergence rates on the physical parameters of viscosity $\mu$ and permeability $\kappa$ can be reduced as well. Similar approaches will be investigated in the three-dimensional case.

In many relevant applications, the porous medium has a heterogeneous structure, and the BIEM is not a suitable choice to model the flow there. This work can be extended to model such porous media, by coupling the BIEM in the domain of free flow with more appropriate techniques such as the finite element method (FEM) in the porous domain with variable properties. With this approach, the FEM can handle the porous domain where the heterogeneities occur, while the BIEM will efficiently deal with the Stokes domain, which can be large or even infinite.

To limit the scope of this paper, only the steady state problem was considered with a single interface. In future work, time-dependent problems and multi-surface cases will be considered. In such a simulation, several surfaces can be close to each other, and it is well-known that the integrands become nearly singular and therefore numerically non-trivial to integrate with sufficient accuracy. The correction methods of~\cite{beale01, beale04, tlupova18} will be used to improve the accuracy in the nearly singular case. Another challenge in such a simulation is the high computational cost of the direct evaluation of the resulting dense linear systems, which is $O(N^2)$, where $N$ is the number of unknowns. To be able to simulate large enough systems, the computational cost will be reduced using fast summation techniques such as the kernel-independent treecode~\cite{wang-krasny-tlupova2} or a fast multipole method~\cite{tornberg-greengard-08}.


\section*{Acknowledgments}
The author gratefully acknowledges support by the National Science Foundation under grant DMS-2012371.

\bibliographystyle{plain}
\bibliography{Bib}

\end{document}